\documentclass[11pt]{article}
\usepackage{amsmath, amssymb, amsfonts, amsthm, amscd, vmargin}
\usepackage[latin1]{inputenc}
\usepackage[all]{xy}
\usepackage{graphicx,color}
\usepackage{enumerate}
\usepackage{tikz}
\usetikzlibrary{shapes}

\setmarginsrb{2.4cm}{3cm}{2.4cm}{3cm}{0cm}{0cm}{0cm}{1.5cm}

\theoremstyle{plain}
\newtheorem{teo}{Theorem}[section]

\newtheorem{prop}[teo]{Proposition}
\newtheorem{cor}[teo]{Corollary}

\theoremstyle{definition}
\newtheorem{defi}[teo]{Definition}
\newtheorem{ex}[teo]{Example}

\newtheorem{rem}[teo]{Remark}

\newtheorem{remdefi}[teo]{Remark/Definition}

\newcommand{\Hom}{\operatorname{Hom}}

\renewcommand{\dim}{\operatorname{dim}}

\renewcommand{\Im}{\operatorname{Im}}

\newcommand{\Cbb}{{\mathbb C}}
\newcommand{\Qbb}{{\mathbb Q}}
\newcommand{\Zbb}{{\mathbb Z}}

\newcommand{\Pbb}{{\mathbb P}}
\newcommand{\Nbb}{{\mathbb N}}

\newcommand{\D}{\Delta}

\def\chadomI{
\draw (0,0) -- (0,10);
\draw[color=gray!60] (1,0) -- (1,10);
\draw[color=gray!60] (2,0) -- (2,10);
\draw[color=gray!60] (3,0) -- (3,10);
\draw[color=gray!60] (4,0) -- (4,10);
\draw[color=gray!60] (5,0) -- (5,10);

\draw[color=gray!60] (0,1) -- (6,1);
\draw[color=gray!60] (0,2) -- (6,2);
\draw[color=gray!60] (0,3) -- (6,3);
\draw[color=gray!60] (0,4) -- (6,4);
\draw[color=gray!60] (0,5) -- (6,5);
\draw[color=gray!60] (0,6) -- (6,6);
\draw[color=gray!60] (0,7) -- (6,7);
\draw[color=gray!60] (0,8) -- (6,8);
\draw[color=gray!60] (0,9) -- (6,9);
\node at (-0.8,2) [color=gray] {$\varpi_0$};
\node at (1.2,0.5) [color=gray] {$\varpi_\alpha$};
\node at (-0.5,1) [color=gray] {0};
\node at (0,1) [color=gray!80] {\textbullet};
\node at (1,1) [color=gray!80] {\textbullet};
\node at (0,2) [color=gray!80] {\textbullet};
}

\def\chadomIbis{
\draw (0,0) -- (0,7);
\draw[color=gray!60] (1,0) -- (1,7);
\draw[color=gray!60] (2,0) -- (2,7);
\draw[color=gray!60] (3,0) -- (3,7);
\draw[color=gray!60] (4,0) -- (4,7);
\draw[color=gray!60] (5,0) -- (5,7);

\draw[color=gray!60] (0,1) -- (6,1);
\draw[color=gray!60] (0,2) -- (6,2);
\draw[color=gray!60] (0,3) -- (6,3);
\draw[color=gray!60] (0,4) -- (6,4);
\draw[color=gray!60] (0,5) -- (6,5);
\draw[color=gray!60] (0,6) -- (6,6);

\node at (-0.8,2) [color=gray] {$\varpi_0$};
\node at (1.2,0.5) [color=gray] {$\varpi_\alpha$};
\node at (-0.5,1) [color=gray] {0};
\node at (0,1) [color=gray!80] {\textbullet};
\node at (1,1) [color=gray!80] {\textbullet};
\node at (0,2) [color=gray!80] {\textbullet};
}

\def\chadomIter{
\draw (0,-2) -- (0,7);
\draw[color=gray!60] (1,-2) -- (1,7);
\draw[color=gray!60] (2,-2) -- (2,7);
\draw[color=gray!60] (3,-2) -- (3,7);
\draw[color=gray!60] (4,-2) -- (4,7);
\draw[color=gray!60] (5,-2) -- (5,7);
\draw[color=gray!60] (0,-1) -- (6,-1);
\draw[color=gray!60] (0,0) -- (6,0);
\draw[color=gray!60] (0,1) -- (6,1);
\draw[color=gray!60] (0,2) -- (6,2);
\draw[color=gray!60] (0,3) -- (6,3);
\draw[color=gray!60] (0,4) -- (6,4);
\draw[color=gray!60] (0,5) -- (6,5);
\draw[color=gray!60] (0,6) -- (6,6);

\node at (-0.8,2) [color=gray] {$\varpi_0$};
\node at (1.2,0.5) [color=gray] {$\varpi_\alpha$};
\node at (-0.5,1) [color=gray] {0};
\node at (0,1) [color=gray!80] {\textbullet};
\node at (1,1) [color=gray!80] {\textbullet};
\node at (0,2) [color=gray!80] {\textbullet};
}

\def\chadomII{
\draw (0,0) -- (0,8);
\draw (0,0) -- (7,0);
\draw[color=gray!80] (0,0) -- (4,8);
\node at (4.7,7.3) [color=gray] {$M_\Qbb$};
\draw[color=gray!60] (1,0) -- (1,8);
\draw[color=gray!60] (2,0) -- (2,8);
\draw[color=gray!60] (3,0) -- (3,8);
\draw[color=gray!60] (4,0) -- (4,8);
\draw[color=gray!60] (5,0) -- (5,8);
\draw[color=gray!60] (6,0) -- (6,8);

\draw[color=gray!60] (0,1) -- (7,1);
\draw[color=gray!60] (0,2) -- (7,2);
\draw[color=gray!60] (0,3) -- (7,3);
\draw[color=gray!60] (0,4) -- (7,4);
\draw[color=gray!60] (0,5) -- (7,5);
\draw[color=gray!60] (0,6) -- (7,6);
\draw[color=gray!60] (0,7) -- (7,7);

\node at (-0.8,1) [color=gray] {$\varpi_\beta$};
\node at (1.2,-0.8) [color=gray] {$\varpi_\alpha$};
\node at (-0.5,-0.5) [color=gray] {0};
\node at (0,0) [color=gray!80] {\textbullet};
\node at (1,0) [color=gray!80] {\textbullet};
\node at (0,1) [color=gray!80] {\textbullet};
}

\def\chadomIIbis{
\draw (0,0) -- (0,6);
\draw (0,0) -- (7,0);
\draw[color=gray!80] (0,0) -- (3,6);
\node at (3.7,5.3) [color=gray] {$M_\Qbb$};
\draw[color=gray!60] (1,0) -- (1,6);
\draw[color=gray!60] (2,0) -- (2,6);
\draw[color=gray!60] (3,0) -- (3,6);
\draw[color=gray!60] (4,0) -- (4,6);
\draw[color=gray!60] (5,0) -- (5,6);
\draw[color=gray!60] (6,0) -- (6,6);

\draw[color=gray!60] (0,1) -- (7,1);
\draw[color=gray!60] (0,2) -- (7,2);
\draw[color=gray!60] (0,3) -- (7,3);
\draw[color=gray!60] (0,4) -- (7,4);
\draw[color=gray!60] (0,5) -- (7,5);

\node at (-0.8,1) [color=gray] {$\varpi_\beta$};
\node at (1.2,-0.8) [color=gray] {$\varpi_\alpha$};
\node at (-0.5,-0.5) [color=gray] {0};
\node at (0,0) [color=gray!80] {\textbullet};
\node at (1,0) [color=gray!80] {\textbullet};
\node at (0,1) [color=gray!80] {\textbullet};
}

\def\chadomIIter{
\draw (0,0) -- (0,4);
\draw (0,0) -- (4,0);
\draw[color=gray!80] (0,0) -- (2,4);
\node at (2.7,3.3) [color=gray] {$M_\Qbb$};
\draw[color=gray!60] (1,0) -- (1,4);
\draw[color=gray!60] (2,0) -- (2,4);
\draw[color=gray!60] (3,0) -- (3,4);

\draw[color=gray!60] (0,1) -- (4,1);
\draw[color=gray!60] (0,2) -- (4,2);
\draw[color=gray!60] (0,3) -- (4,3);

\node at (-0.5,1) [color=gray] {$\varpi_\beta$};
\node at (1.2,-0.5) [color=gray] {$\varpi_\alpha$};
\node at (-0.3,-0.3) [color=gray] {0};
\node at (0,0) [color=gray!80] {\textbullet};
\node at (1,0) [color=gray!80] {\textbullet};
\node at (0,1) [color=gray!80] {\textbullet};
}

\begin{document}

\title{The Log Minimal Model Program for horospherical varieties via moment polytopes}

\author{Boris Pasquier\protect\footnote{Boris PASQUIER,
Institut Montpelli\'erain Alexander Grothendieck, CNRS, Univ. Montpellier.
E-mail: boris.pasquier@umontpellier.fr}}

\maketitle
\begin{abstract}
In \cite{MMPhoro}, we described the Minimal Model Program in the family of $\Qbb$-Gorenstein projective horospherical varieties, by studying certain continuous changes of moment polytopes of polarized horospherical varieties.  Here, we summarize the results of \cite{MMPhoro} and we explain how to generalize them in order to describe the Log Minimal Model Program for pairs $(X,\D)$ when $X$ is a projective horospherical variety. 
\end{abstract}



	\section{Introduction}
	
	In this paper, we work over the complex numbers. And all varieties are supposed to be irreducible.
	
	Let $G$ be a connected reductive algebraic group. Let $X$ be a normal $G$-variety. We say that $X$ is horospherical if there exists $x\in X$ such that $G\cdot x$ is open in $X$ and $x$ is fixed by a maximal unipotent subgroup of $G$. 

Note that if $X$ is horospherical then there exists a Borel subgroup $B$ of $G$ such that $B\cdot x$ is open in $X$, i.e., $X$ is a spherical variety.  Moreover, if $X$ is projective and horospherical, then $X$ is the closure of the $G$-orbit of a sum of highest weight vectors in the projectivization of a multiplicity free $G$-module (see \cite[Proposition~2.11]{MMPhoro}). This point of view motivates the classification of projective horospherical varieties in terms of polytopes (called moment polytopes). In this paper, we propose to describe different variations of the Minimal Model Program (MMP) for projective horospherical varieties (including projective toric varieties) via moment polytopes.\\

The principle of the MMP is to construct representatives in birational equivalence classes of projective varieties, called minimal models. In fact, a minimal model does not always exist, for example for rational varieties;  then the MMP constructs either minimal models or certain types of fibrations called Mori fibrations. The  idea is to contract ``negative curves'', and if neccesary (i.e., if singularities are too ``bad'') to partially  desingularize by adding ``positive'' curves, in order to terminate either with a minimal model, which has only ``non-negative'' curves, or with a Mori fibration, whose general fiber has only ``negative'' curves.\\
 The ``sign'' of a curve, in the classical MMP, is given by the sign of its intersection with a canonical divisor of the variety. This intersection is well-defined only when the canonical divisor is $\Qbb$-Cartier (i.e., the variety is $\Qbb$-Gorenstein). The Log MMP consists on adding a (``small'') $\Qbb$-divisor to the canonical divisor, so that the sum is $\Qbb$-Cartier.

We now describe the MMP and the Log MMP in more details.

\begin{defi}\label{def:pairs} Let $X$ be a normal projective variety. Let $\D$ be a $\Qbb$-divisor (not necessarily $\Qbb$-Cartier).
Let $K_X$ be a canonical divisor of $X$. We say that $(X,\D)$ is a \emph{log pair} if $K_X+\D$ is $\Qbb$-Cartier.
\end{defi}
Note that the divisor $\D$ is sometimes supposed to be effective in the Log MMP theory. Here, we do not make this assumption.\\

We denote by $NE(X)$ the cone of numerical classes of effective curves on $X$, by $\overline{NE}(X)$ its closure and by $\overline{NE}(X)_{K_X+\D<0}$ (respectively $\overline{NE}(X)_{K_X+\D>0}$) the open half-space of $\overline{NE}(X)$ defined by curves that are negative (respectively positive) along the divisor $K_X+\D$.
 
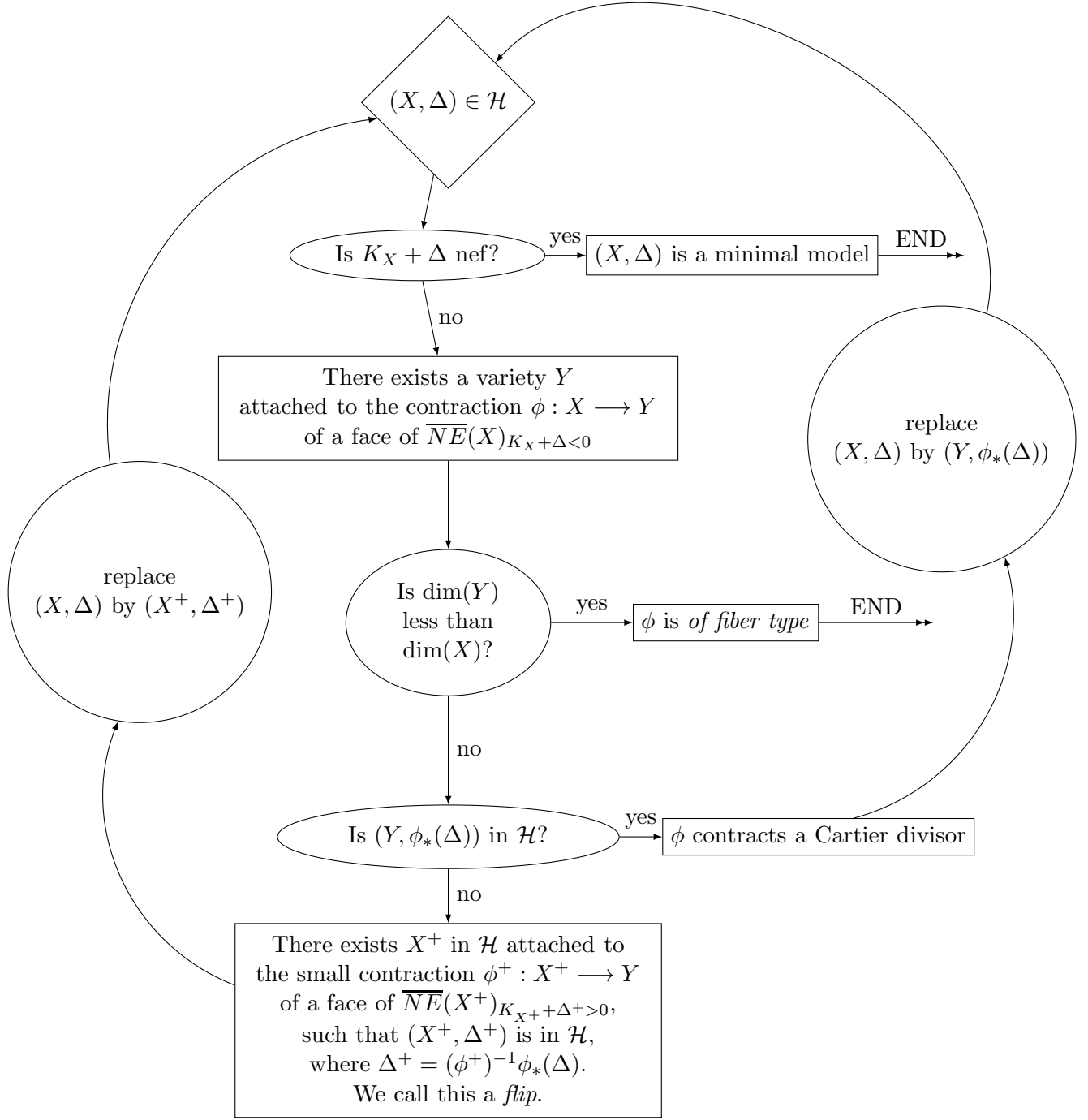
\begin{figure}
\caption{Log MMP in a family $\mathcal{H}$ of log pairs}
\label{fig:logMMP}
\begin{tikzpicture}[scale=1]
\node [draw][diamond] (a) at (0,0) {
   $(X,\D)\in\mathcal{H}$
};
\node [draw][ellipse] (b) at (-0.5,-2.5) {
   Is $K_X+\D$ nef?
};
\node [draw][rectangle] (c) at (4.6,-2.5) {
   $(X,\D)$ is a minimal model
};
\node [rectangle] (d) at (8.5,-2.5) {};
\node [draw] [rectangle] (e) at (0,-5) {
\begin{tabular}{c}
There exists a variety $Y$\\ attached to the contraction $\phi:X\longrightarrow Y$\\ of a face of $\overline{NE}(X)_{K_X+\D<0}$
\end{tabular}
};
\node [draw] [ellipse] (f) at (0,-8.5) {
\begin{tabular}{c}
Is $\dim(Y)$\\ less than \\ $\dim(X)$?
\end{tabular}
};
\node [draw] [rectangle] (g) at (4.5,-8.5) {
$\phi$ is \emph{of fiber type}
};
\node [rectangle] (h) at (8,-8.5) {};
\node [draw] [ellipse] (i) at (0,-12) {
\begin{tabular}{c}
Is $(Y,\phi_*(\D))$ in $\mathcal{H}$?
\end{tabular}
};
\node [draw] [rectangle] (j) at (6,-12) {
$\phi$ contracts a Cartier divisor
};
\node [draw] [rectangle] (k) at (0,-15) {
\begin{tabular}{c}
There exists $X^+$ in $\mathcal{H}$ attached to\\ the small contraction $\phi^+:X^+\longrightarrow Y$\\ of a face of $\overline{NE}(X^+)_{K_{X^+}+\D^+>0}$,\\
such that $(X^+,\D^+)$ is in $\mathcal{H}$,\\
where $\D^+=(\phi^+)^{-1}\phi_*(\D)$.\\
We call this a \emph{flip}.
\end{tabular}
};
\node [draw][circle] (l) at (8,-5.5) {\begin{tabular}{c} replace \\ $(X,\D)$ by $(Y,\phi_*(\D))$\end{tabular}};
\node [draw][circle] (m) at (-5,-8) {\begin{tabular}{c} replace \\ $(X,\D)$ by $(X^+,\D^+)$\end{tabular}};

\draw[->,>=latex] (a) -- (b);
\draw[->,>=latex] (b) -- node[above] {yes} (c);
\draw[->>,>=latex] (c) --  node[above] {END} (d);
\draw[->,>=latex] (b) --  node[right] {no} (e);
\draw[->,>=latex] (e) -- (f);
\draw[->,>=latex] (f) --  node[above] {yes} (g);
\draw[->>,>=latex] (g) --  node[above] {END} (h);
\draw[->,>=latex] (f) --  node[right] {no} (i);
\draw[->,>=latex] (i) --  node[above] {yes} (j);
\draw[->,>=latex] (i) --  node[right] {no} (k);

\draw[->,>=latex] (j) to[bend right=45] (l);
\draw[->,>=latex] (l) to[bend right=75] (a);
\draw[->,>=latex] (k) to[bend left=45] (m);
\draw[->,>=latex] (m) to[bend left=45] (a);

\end{tikzpicture}
\end{figure}

We summarized, in Figure~\ref{fig:logMMP}, the principle of the Log MMP (without $\Qbb$-factorial assumption).

The different versions of MMP and Log MMP depend on the choice of the family  $\mathcal{H}$:
\begin{itemize}
\item $\Qbb$-factorial MMP (or generally just called MMP), when $\mathcal{H}$ is the set of log pairs $(X,0)$ such that $X$ is $\Qbb$-factorial.
\item non-$\Qbb$-factorial MMP (or $\Qbb$-Gorenstein MMP), when $\mathcal{H}$ is the set of log pairs $(X,0)$ such that $X$ is $\Qbb$-Gorenstein.
\item $\Qbb$-factorial Log MMP (or generally just called Log MMP), when $\mathcal{H}$ is the set of log pairs $(X,\D)$ such that $X$ is $\Qbb$-factorial.
\item non-$\Qbb$-factorial Log MMP, when $\mathcal{H}$ is the set of all log pairs $(X,\D)$.
\end{itemize}

In \cite{MMPhoro}, the first two families were considered for horospherical varieties by reducing to the description of one-parameter families of polytopes. Moreover, in \cite{MMPGor}, the non-$\Qbb$-factorial MMP and Log MMP were also detailed in general before to discuss the case of spherical varieties. 

In this paper, we consider the last two families when $X$ is horospherical.

Note that, in these two families, we can distinguish several subfamilies according to types of singularities: terminal, canonical,  Kawamata log terminal (klt), purely log terminal, divisorial log terminal, weakly log-terminal and log canonical (lc). For more details about these different types of singularities see for example \cite{Fujino4}. Here, we will only deal with klt and log canonical singularities, whose definitions (see Definition~\ref{def:kltlc}) do not depend on the log resolution (by assuming that the coefficients of $\D$ are at most one in the definition of lc singularities). 

The results that we obtain in this paper can be summarized in the following theorem.

\begin{defi}\label{def:horopair}
A \emph{horospherical pair} is a log pair $(X,\D)$ such that $X$ is a horospherical $G$-variety and $\D$ is $B$-stable, where $G$ and $B$ are chosen as explained above. 

\end{defi}

\begin{teo}\label{th:main0}
Let $(X,\D)$ be a projective horospherical pair. Suppose that $(K_X+\D)$ is non-zero.

Then, for any choice of an ample $\Qbb$-Cartier $\Qbb$-divisor of $X$, we can construct a family of polytopes, depending on one non-negative rational parameter, such that: each polytope is either the empty set or the moment polytope of a projective horospherical variety; there are finitely many horospherical varieties that correspond to a polytope of the family; and all these horospherical varieties are the ones that appear in the Log MMP.

In particular the Log MMP described by Figure~\ref{fig:logMMP} works and terminates as soon as $-(K_X+\D)$ is non-zero,  if $\mathcal{H}$ is the family of klt horospherical pairs, if $\mathcal{H}$ is the family of lc horospherical pairs, or if $\mathcal{H}$ is the family of horospherical pairs $(X,\D)$ such that $(K_X+\D)$ is effective. We can also restrict these families to log pairs with $\D$ effective.

Moreover,  we can explicitly describe each step of the Log MMP until it terminates.
\end{teo}

The strategy is the same as in \cite{MMPhoro}  and the proof of Theorem~\ref{th:main0} is very similar to the proof of \cite[Theorem~1.1]{MMPhoro}.  Hence, in Section~\ref{sec:horopoly}, we follow the notation of \cite{MMPhoro}  and we recall briefly the definitions and properties that we need here. Then, in Section~\ref{sec:FamilyConstruction}, we explain why the proofs of \cite{MMPhoro} can be easily adapted to the log case. And we conclude and give more examples in Section~\ref{sec:conclusion}.

\section{Projective horospherical varieties and polytopes}\label{sec:horopoly}

\subsection{Notation}

We will not recall here the long Luna-Vust theory of horospherical embeddings (classification in terms of colored fans).  For more details on horospherical varieties, we refer the reader to \cite{Fanohoro}, and for basic results on Luna-Vust theory of spherical embeddings, we refer to \cite{Knop89}. 

In this paper, we will only describe and use another classification of projective horospherical varieties in terms of moment polytopes, which was first introduced in \cite[Section 2.3]{MMPhoro}.\\

We fix a connected reductive algebraic group $G$, a Borel subgroup $B$ of $G$ and a maximal torus $T$  in $B$.

Then $G/H$ always denotes a horospherical homogeneous space such that $H$ contains the unipotent radical of $B$. And if $X$ is a horospherical variety, we denote by $G/H$ the horospherical homogeneous space as above that is isomorphic to the open $G$-orbit of $X$, and we say that $X$ is a $G/H$-embedding. (Note that Borel subgroups are all conjugated in $G$, so that the assumption on $H$ can be done without loss of generality.) 

In order to simplify, we say here that a \emph{$G/H$-embedding} is a normal algebraic $G$-variety with an open $G$-orbit isomorphic to $G/H$. (See \cite[section 2.1]{MMPhoro} for the precise definition and the notion of isomorphism of $G/H$-embeddings.)\\

The normalizer of $H$ in $G$ is a parabolic subgroup of $G$, we denote it by $P$.

Let $S$ be the set of simple roots of $(G,B,T)$. Then, we denote by $R$ the subset of $S$ of simple roots of $P$.
Let $X(T)$ (respectively $X(T)^+$) be the lattice of characters of $T$ (respectively the set of dominant characters). Similarly, we define $X(P)$ and $X(P)^+=X(P)\cap X(T)^+$. Note that the lattice $X(P)$ and the cone $X(P)^+$ are generated by the fundamental weights $\varpi_\alpha$ with $\alpha\in S\backslash R$ and the weights of the center of $G$.\\

We denote by $M$ the sublattice of $X(P)$ consisting of characters of $P$ vanishing on $H$. The rank of $M$ is called the \emph{rank} of $G/H$, and denoted by $n$. Let $N:=\Hom_\Zbb(M,\Zbb)$ be the dual of~$M$.

For any free lattice $\mathbb{L}$, we denote by $\mathbb{L}_\Qbb$ the $\Qbb$-vector space $\mathbb{L}\otimes_\Zbb\Qbb$.\\

 For any $\alpha\in S\backslash R$, we define $$W_{\alpha,P}:=\{m\in X(P)_\Qbb\,\mid\,\langle m,\alpha^\vee\rangle=0\}.$$ Note that these hyperplanes correspond to the walls of the dominant chamber $X(P)^+$.

Throughout the paper, we illustrate the results with two main examples of horospherical homogeneous spaces.

\begin{ex}\label{ex:GH}
~
\begin{itemize}
\item {\it\large Example~I}:
let $G=\rm{SL}_2(\Cbb)\times\Cbb^*$ and let $H$ be the unipotent radical  $\lbrace\left(\left(\begin{array}{cc}1 & x\\0 & 1
\end{array}\right),1\right)\,\mid\,x\in\Cbb\rbrace$ of the Borel subgroup of $G$ defined by $B=\lbrace\left(\left(\begin{array}{cc}t & x\\0 & t^{-1}
\end{array}\right),u\right)\,\mid\,x\in\Cbb,\,\,t,u\in\Cbb^*\rbrace.$ Then $P=B$ and let $M=X(P)$. We denote by $\varpi_\alpha$ the fundamental weight of $(G,B)$ associated to the unique simple root of $(G,B)$. And we denote by $\varpi_0$ the weight of $B$ defined by the projection on $\Cbb^*$. Then $M=\Zbb\varpi_\alpha\oplus\Zbb\varpi_0$ and $X(P)^+=\Nbb\varpi_\alpha\oplus\Zbb\varpi_0$. Note also that $n=2$ and $G/H$ is three-dimensional.
\item {\it\large Example~II}:
let $G=\rm{SL}_3(\Cbb)$ and $H=\lbrace\left(\begin{array}{ccc}t^2 & x & y\\0 & t^{-3} & z\\0& 0 & t
\end{array}\right)\,\mid\,x,y,z\in\Cbb,\,\,t\in\Cbb^*\rbrace.$ Then $P$ is the Borel subgroup $B$ of $G$ consisting of upper triangular matrices and $M=\Zbb(\varpi_\alpha+2\varpi_\beta)$, where $\varpi_\alpha$ and $\varpi_\beta$ are the two fundamental weights of $(G,B)$ (defined by $\varpi_\alpha(b)=b_{11}$ and $\varpi_\beta(b)=b_{11}b_{22}$ for $b\in B$). Moreover, $X(P)=\Zbb\varpi_\alpha\oplus\Zbb\varpi_\beta$ and $X(P)^+=\Nbb\varpi_\alpha\oplus\Nbb\varpi_\beta$. Note also that $n=1$, $G/H$ is four-dimensional and $M_\Qbb\neq X(P)_\Qbb$.
\end{itemize}
\end{ex}

\subsection{Divisors of horospherical varieties and moment polytopes}\label{sec:divpoly}

In this section, we explain how to construct a polytope from an ample $\Qbb$-Cartier $\Qbb$-divisor of a projective horospherical variety.

First, let us justify why we only consider $B$-stable divisors.

\begin{prop}(\cite[Section 2.2]{briondiv}) \label{prop:divBstable}
Any divisor of any $G/H$-embedding is linearly equivalent to a $B$-stable divisor.
\end{prop}

Now, we describe the $B$-stable prime divisors of a $G/H$-embedding $X$. We denote by $X_1,\dots,X_r$ the $G$-stable prime divisors of $X$. The other $B$-stable prime divisors (i.e., those that are not $G$-stable) are the closures in $X$ of $B$-stable prime divisors of $G/H$ (which are called \emph{colors} of $G/H$). They are indexed by the simple roots $\alpha$ in $S\backslash R$ and we denote them by $D_\alpha$ (see \cite[Section~2]{Fanohoro} for an explicit description of these divisors).

Hence, any $B$-stable divisor (respectively $\Qbb$-divisor) of $X$ can be written as follows:\\ $\sum_{i=1}^r d_i X_i +\sum_{\alpha\in S\backslash R} d_\alpha D_\alpha$ with the $d_i$'s and the $d_\alpha$'s in $\Zbb$ (respectively $\Qbb$).\\

All these $B$-stable prime divisors have an image in the lattice $N$ as follows. Note that the lattice $M$ is in bijection with the set of $B$-semi-invariant rational functions of $G/H$ up to a scalar multiple. Any $B$-stable divisor induces a $B$-stable valuation and then, by restriction to $B$-semi-invariant rational functions, it induces a point in $N$. We denote by $x_i$ the corresponding point of $N$ associated to $X_i$, it is a non-zero primitive element of the lattice $N$.
For any simple root $\alpha\in S\backslash R$, the point associated to $D_\alpha$ is the restriction (maybe zero) of the coroot $\alpha^\vee$ to $M$, which we denote by $\alpha^\vee_M$ (see \cite[Section~2]{Fanohoro}).\\

In \cite{briondiv}, there are characterizations of Cartier, $\Qbb$-Cartier and ample $B$-stable divisors of spherical varieties and also a description of the global sections of such divisors in terms of polytopes, that permits to give the following definition and result.

In the rest of the section, we fix a projective $G/H$-embedding $X$ and a $\Qbb$-divisor $D=\sum_{i=1}^r d_i X_i +\sum_{\alpha\in S\backslash R} d_\alpha D_\alpha$ of $X$. And we suppose that $D$ is $\Qbb$-Cartier and ample. 
\begin{defi}
The \emph{pseudo-moment polytope} of $(X,D)$ is $$\tilde{Q}_D:=\{m\in M_\Qbb\,\mid\,\langle m,x_i\rangle\geq -d_i,\,\forall i\in\{1,\dots,r\}\mbox{ and }\langle m,\alpha^\vee_M\rangle\geq -d_\alpha,\,\forall \alpha\in S\backslash R \}.$$
Let $v^0:=\sum_{\alpha\in S\backslash R}d_\alpha\varpi_\alpha$. The \emph{moment polytope}  of $(X,D)$ is $Q_D:=v^0+\tilde{Q}_D$.
\end{defi}

Note that $v^0$ is not necessarily in $M_\Qbb$, but only in $X(P)_\Qbb$.

\begin{prop}\cite[Corollary 2.8]{MMPhoro}\label{prop:moment}
\begin{enumerate}
\item The pseudo-moment polytope $\tilde{Q}_D$ of $(X,D)$ is of maximal dimension in $M_\Qbb$.
\item The moment polytope $Q_D$ is contained in the dominant chamber $X(P)^+_\Qbb$ and it is not contained in any wall $W_{\alpha,P}$ for $\alpha\in S\backslash R$. 
\item There is a bijection between faces of $Q_D$ (or $\tilde{Q}_D$) and $G$-orbits of $X$ (preserving the respective orders). In particular, the $G$-stable primes divisors $X_i$ are in bijection with the facets of $Q_D$ that are not contained in any wall $W_{\alpha,P}$ for $\alpha\in S\backslash R$ (the bijection maps $X_i$ to the facet of $\tilde{Q}_D$ defined by $\langle m,x_i\rangle =-d_i$).
\item The divisor $D$ can be computed from the pair $(Q,\tilde{Q})$ as follows: the coefficients $d_\alpha$ with $\alpha\in S\backslash R$ are given by the translation vector in $X(P)_\Qbb$ that maps $\tilde{Q}$ to $Q$; and for any $i\in\{1,\dots,r\}$, the coefficient $d_i$ is given by $-\langle v_i,x_i\rangle$ for any element $v_i\in M_\Qbb$ in the facet of $\tilde{Q}$ associated to $X_i$.
\end{enumerate}
\end{prop}

\subsection{Classification of projective horospherical varieties in terms of polytopes}\label{sec:momentpoly}

In this section, we write the classification of projective $G/H$-embeddings in terms of $G/H$-polytopes (defined below in Definition~\ref{def:G/H-eq}) and we also give a similar classification of polarized projective horospherical varieties.


\begin{defi}\label{def:G/H-eq}
Let $Q$ be a polytope in $X(P)^+_\Qbb$ (not necessarily a lattice polytope). We say that $Q$ is a \emph{$G/H$-polytope}, if its direction is $M_\Qbb$ (i.e., if the subvector space of  $X(P)_\Qbb$, spanned by the vectors $x-y$ with $x$ and $y$ in $Q$, is $M_\Qbb$) and if it is contained in no wall $W_{\alpha,P}$ with $\alpha\in S\backslash R$.

Let $Q$ and $Q'$ be two $G/H$-polytopes in $X(P)^+_\Qbb$. Consider any polytopes $\tilde{Q}$ and $\tilde{Q'}$ in $M_\Qbb$ obtained by translations from $Q$ and $Q'$ respectively.   We say that $Q$ and $Q'$ are \emph{equivalent $G/H$-polytopes} if the following conditions are satisfied.
\begin{enumerate}
\item There exist an integer $j$ and $2j$ affine half-spaces $\mathcal{H}_1^+,\dots,\mathcal{H}_j^+$ and $\mathcal{H'}_1^+,\dots,\mathcal{H'}_j^+$ of $M_\Qbb$ (respectively delimited by the affine hyperplanes $\mathcal{H}_1,\dots,\mathcal{H}_j$ and $\mathcal{H'}_1,\dots,\mathcal{H'}_j$) such that $\tilde{Q}$ is the intersection of the $\mathcal{H}_i^+$, $\tilde{Q'}$ is the intersection of the $\mathcal{H'}_i^+$, and for all $i\in\{1,\dots,j\}$, $\mathcal{H}_i^+$ is the image of $\mathcal{H'}_i^+$ by a translation.
\item With the notation of the previous item, for all subsets $J$ of $\{1,\dots,j\}$, the intersections $\cap_{i\in J}\mathcal{H}_i\cap \tilde{Q}$ and $\cap_{i\in J}\mathcal{H'}_i\cap \tilde{Q'}$ have the same dimension.
\item $Q$ and $Q'$ intersect exactly the same walls $W_{\alpha,P}$ of $X(P)^+_\Qbb$ (with $\alpha\in S\backslash R$).
\end{enumerate}
\end{defi}

\begin{prop}\cite[Proposition~2.10]{MMPhoro}\label{prop:G/H-eq}
The map from (isomorphism classes of) projective $G/H$-embeddings to the set of equivalence classes of $G/H$-polytopes that maps $X$ to the class of the moment polytope of $(X,D)$, where $D$ is any ample $\Qbb$-Cartier $B$-stable $\Qbb$-divisor, is a well-defined bijection.
\end{prop}

\begin{ex}\label{ex:Qi}
We consider the homogeneous spaces $G/H$ of Example~\ref{ex:GH}.
\begin{itemize}
\item {\it\large Example~I}: 
the six $G/H$-polytopes below correspond to four disctinct projective $G/H$-embeddings. The $G/H$-polytopes $Q_1$ and $Q_2$ are equivalent, as well as $Q_4$ and $Q_5$. 
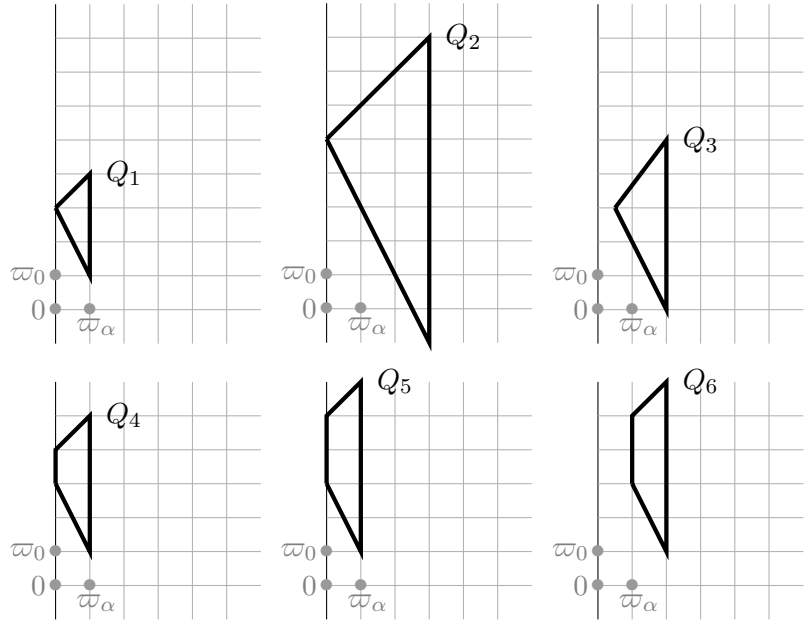
\begin{center}

\begin{tikzpicture}[scale=0.45]
\chadomI
\draw[ultra thick] (0,4) -- (1,5) -- (1,2) -- (0,4);
\node at (2,5) {$Q_1$};
\end{tikzpicture}
\begin{tikzpicture}[scale=0.45]
\chadomI
\draw[ultra thick] (0,6) -- (3,9) -- (3,0) -- (0,6);
\node at (4,9) {$Q_2$};
\end{tikzpicture}
\begin{tikzpicture}[scale=0.45]
\chadomI
\draw[ultra thick] (1/2,4) -- (2,6) -- (2,1) -- (1/2,4);
\node at (3,6) {$Q_3$};
\end{tikzpicture}

\begin{tikzpicture}[scale=0.45]
\chadomIbis
\draw[ultra thick] (0,4) -- (0,5) -- (1,6) -- (1,2) -- (0,4);
\node at (2,6) {$Q_4$};
\end{tikzpicture}
\begin{tikzpicture}[scale=0.45]
\chadomIbis
\draw[ultra thick] (0,4) -- (0,6) -- (1,7) -- (1,2) -- (0,4);
\node at (2,7) {$Q_5$};
\end{tikzpicture}
\begin{tikzpicture}[scale=0.45]
\chadomIbis
\draw[ultra thick] (1,4) -- (1,6) -- (2,7) -- (2,2) -- (1,4);
\node at (3,7) {$Q_6$};
\end{tikzpicture}
\end{center}

\item {\it\large Example~II}:
the six $G/H$-polytopes below correspond to four disctinct projective $G/H$-embeddings. The $G/H$-polytopes $Q_7$ and $Q_8$ are equivalent, as well as $Q_9$ and $Q_{10}$. 
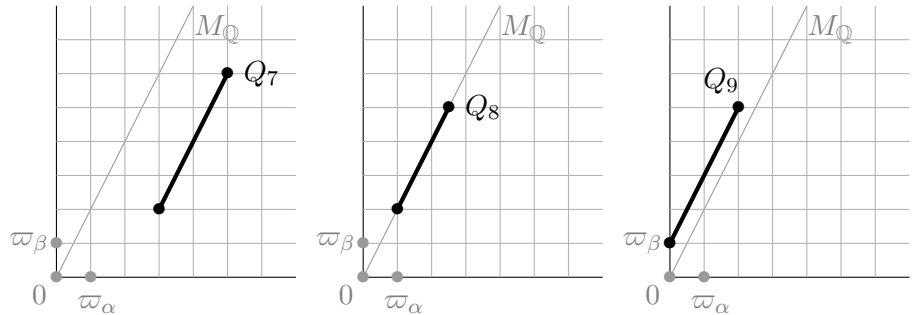
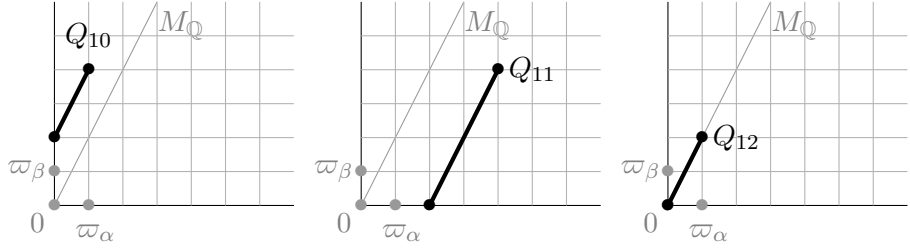
\begin{center}
\begin{tikzpicture}[scale=0.45]
\chadomII
\draw[ultra thick] (5,6) -- (3,2);
\node at (5,6)  {\textbullet};
\node at (3,2)  {\textbullet};
\node at (6,6) {$Q_7$};
\end{tikzpicture}
\begin{tikzpicture}[scale=0.45]
\chadomII
\draw[ultra thick] (5/2,5) -- (1,2);
\node at (5/2,5) {\textbullet};
\node at (1,2) {\textbullet};
\node at (7/2,5) {$Q_8$};
\end{tikzpicture}
\begin{tikzpicture}[scale=0.45]
\chadomII
\draw[ultra thick] (0,1) -- (2,5);
\node at (0,1) {\textbullet};
\node at (2,5) {\textbullet};
\node at (1.5,5.8) {$Q_9$};
\end{tikzpicture}

\begin{tikzpicture}[scale=0.45]
\chadomIIbis
\draw[ultra thick] (0,2) -- (1,4);
\node at (0,2) {\textbullet};
\node at (1,4) {\textbullet};
\node at (1,5) {$Q_{10}$};
\end{tikzpicture}
\begin{tikzpicture}[scale=0.45]
\chadomIIbis
\draw[ultra thick] (2,0) -- (4,4);
\node at (2,0) {\textbullet};
\node at (4,4) {\textbullet};
\node at (5,4) {$Q_{11}$};
\end{tikzpicture}
\begin{tikzpicture}[scale=0.45]
\chadomIIbis
\draw[ultra thick] (0,0) -- (1,2);
\node at (0,0) {\textbullet};
\node at (1,2) {\textbullet};
\node at (2,2) {$Q_{12}$};
\end{tikzpicture}
\end{center}
\end{itemize}
\end{ex}

Since isomorphism classes of horospherical homogeneous $G$-spaces are in bijection with pairs $(P,M)$ where $P$ is a parabolic subgroup of $G$ containing $B$ and $M$ is a sublattice of $X(P)$ (see \cite[Proposition~2.4]{MMPhoro}), we can give the following alternative classification.

\begin{defi}\label{def:momenttriple}
A \emph{moment quadruple} is a quadruple $(P,M,Q,\tilde{Q})$ where $P$ is a parabolic subgroup of $G$ containing $B$, $M$ is a sublattice of $X(P)$, $Q$ is a polytope in $X(P)^+_\Qbb$ and $\tilde{Q}$ is a polytope in $M_\Qbb$ that satisfies the three following conditions. 
\begin{enumerate} 
\item There exists (a unique) $\varpi\in\mathfrak{X}(P)_\Qbb$ such that $Q=\varpi+\tilde{Q}$.
\item The polytope $\tilde{Q}$ is of maximal dimension in $M_\Qbb$ (i.e., its interior in $M_\Qbb$ is not empty).\item The polytope $Q$ is not contained in any wall $W_{\alpha,P}$ for $\alpha\in S\backslash R$.
\end{enumerate}
\end{defi}

\begin{defi}
A \emph{polarized horospherical variety} (respectively $G/H$-embedding) is a pair $(X,D)$ such that $X$ is a projective horospherical variety (respectively $G/H$-embedding) and $D$ is an ample $\Qbb$-Cartier $B$-stable $\Qbb$-divisor. We say that $(X,D)$ is isomorphic to $(X',D')$ if there is an isomorphism from $X$ to $X'$ of embeddings of the same homogeneous space such that $D$ gets identified with $D'$ under this isomorphism.
\end{defi}

\begin{cor}\label{cor:polytopepolarized}
The map from the set of isomorphism classes of polarized projective horospherical $G$-varieties to the set of classes of moment quadruples, that maps $(X,D)$ to $(P,M,Q_D,\tilde{Q}_D)$ is a bijection.
\end{cor}
 
\subsection{$G$-equivariant morphisms between projective horospherical varieties and polytopes}\label{sec:morph}

The existence of dominant $G$-equivariant morphisms between projective horospherical varieties can be characterized in terms of colored fans \cite{Knop89} but also in terms of moment polytopes.

Consider two horospherical homogeneous spaces $G/H$ and $G/H'$.
Let $(X,D)$ be a polarized $G/H$-embedding, let $(X',D')$ be a polarized $G/H'$-embedding.   By Corollary~\ref{cor:polytopepolarized}, $(X,D)$ corresponds to a moment quadruple $(P,M,Q,\tilde{Q})$ and $(X',D')$ corresponds to a moment quadruple $(P',M',Q',\tilde{Q'})$. (We also denote by $R'$ the set of simple roots of $P'$ and by $N'$ the dual lattice of $M'$.)\\

A first necessary condition for the existence of a dominant $G$-equivariant morphism from $X$ to $X'$, is the existence of a projection $\pi$ from $G/H$ to $G/H'$. In particular $H'\supset H$, $P'\supset P$ and $R'\supset R$. The projection $\pi$ induces an injective morphism $\pi^*$ from $M'$ to $M$. We suppose that this necessary condition is satisfied in the rest of the section and we identify $M'$ with $\pi^*(M')$.\\

We first need to define a map $\psi$ from the set of facets of $\tilde{Q}$ to the set of faces of $\tilde{Q'}$ (including $\tilde{Q'}$ itself).

\begin{remdefi}
Note a general fact on polytopes: if $\mathcal{P}$ is a polytope in $\Qbb^r$, then for any affine half-space $\mathcal{H}^+$ delimited by an affine hyperplane $\mathcal{H}$ in $\Qbb^r$, there exists a unique face $F$ of $\mathcal{P}$ and a point $x\in\Qbb^r$ such that $F$ is defined by $x+\mathcal{H}$ (i.e., $F=\mathcal{P}\cap (x+\mathcal{H})$ and $\mathcal{P}\subset x+ \mathcal{H}^+$).
Then, for any facet $F$ of $\tilde{Q}$, let  $\mathcal{H}^+$ be the affine half-space in $M_\Qbb$, delimited   by an affine hyperplane $\mathcal{H}$, such that $F=\mathcal{H}\cap\tilde{Q}$ and $\tilde{Q}\subset \mathcal{H}^+$. If $\mathcal{H}^+\cap M'_\Qbb\neq M'_\Qbb$, it is an affine half-space in $M'_\Qbb$ and, applying the above fact to $\mathcal{P}=\tilde{Q'}$ yields a unique face $F'$ of $\tilde{Q'}$. In that case, we define $\psi(F)$ to be $F'$. And if $\mathcal{H}^+\cap M'_\Qbb=M'_\Qbb$,  we set $\psi(F)=\tilde{Q'}$.
\end{remdefi}

\begin{ex}
To illustrate the definition of $\psi$, we consider the poytopes $Q=Q_4$ and $Q'=Q_1$ of Example~\ref{ex:Qi}.
Then, we have $\psi([AB])=[A'B']$, $\psi([BC])=[B'C']$, $\psi([CD])=[A'C']$, and $\psi([AD])=A'$.
\begin{center}
\begin{tikzpicture}[scale=0.45]
\chadomIbis
\draw[ultra thick] (0,4) -- (0,5) -- (1,6) -- (1,2) -- (0,4);
\node at (3,4) {$Q$};
\node at (-0.5,5.2) {$A$};
\node at (1.5,6.3) {$B$};
\node at (1.5,1.7) {$C$};
\node at (-0.5,3.8) {$D$};
\end{tikzpicture}
\begin{tikzpicture}[scale=0.45]
\chadomIbis
\draw[ultra thick] (0,4) -- (1,5) -- (1,2) -- (0,4);
\node at (3,4) {$Q'$};
\node at (-0.5,4.2) {$A'$};
\node at (1.5,5.3) {$B'$};
\node at (1.5,1.7) {$C'$};
\end{tikzpicture}
\end{center}
\end{ex}

See also \cite[Example~2.15]{MMPhoro} to get an example with $Q$ and $Q'$ of different dimensions.  \\

We can now characterize the existence of dominant $G$-equivariant morphisms.

\begin{prop}\cite[Proposition~2.16 and Corollary~2.17]{MMPhoro}\label{prop:morph}
Under the above necessary condition, there exists a dominant $G$-equivariant morphism from $X$ to $X'$, if and only if 
\begin{enumerate}
\item for any set $\mathcal{G}$ of facets of $\tilde{Q}$ with $\cap_{F\in\mathcal{G}}F\neq\emptyset$ we have  $\cap_{F\in\mathcal{G}}\psi(F)\neq\emptyset$, and
\item for any $\alpha\in S\backslash R$ such that $Q\cap W_{\alpha,P}\neq\emptyset$ , we have $Q'\cap W_{\alpha,P}\neq\emptyset$. 
\end{enumerate}

Suppose there exists a dominant $G$-equivariant morphism $\phi$ from $X$ to $X'$. Let $\mathcal{O}$ be the $G$-orbit in $X$ associated to the face $\cap_{F\in\mathcal{G}}F$ for $\mathcal{G}$ a set of facets. Then $\phi(\mathcal{O})$ is the $G$-orbit in $X'$ associated to the face $\cap_{F\in\mathcal{G}}\psi(F)$ of $Q'$.
\end{prop}

Remark that, in Proposition~\ref{prop:morph}~(1), we could replace the pseudo-moment polytope by the moment polytope  (by extending the definition of $\psi$).

\begin{ex}\label{ex:morph}
We consider the $G/H$-polytopes of Example~\ref{ex:Qi}. Denote by $X_i$ the $G/H$-embedding corresponding to $Q_i$. (Recall that $X_1\simeq X_2$, $X_4\simeq X_5$, $X_7\simeq X_8$ and $X_9\simeq X_{10}$.)
\begin{itemize}
\item {\it\large Example~I}: there exist $G$-equivariant morphisms from $X_3$ to $X_1$, from $X_4$ to $X_1$, and from $X_6$ to $X_1$, $X_3$ and $X_4$.

\item {\it\large Example~II}: there exist $G$-equivariant morphisms from $X_7$ to $X_9$, $X_{11}$ and $X_{12}$, from $X_9$ to $X_{12}$, and from $X_{11}$ to $X_{12}$.

\end{itemize}
\end{ex}

\subsection{Curves in horospherical varieties}\label{sec:curves}

In \cite{Brionmori}, M.~Brion describes the curves of spherical varieties and their intersections with divisors. We rewrite here some of his results for horospherical varieties and in terms of moment polytopes. In particular, we describe effective curves of a projective horospherical variety $X$ and their intersections with ample $\Qbb$-Cartier $\Qbb$-divisors. 
 
We denote by $N_1(X)$ the group of numerical classes of 1-cycles of the variety $X$. Recall that $NE(X)$ is the convex cone in $N_1(X)$ generated by effective 1-cycles.\\

\begin{prop}\cite[Section~2.5]{MMPhoro}
Let $X$ be a projective horospherical variety.
Let $Q$ be any moment polytope of $X$ (with any choice of an ample $\Qbb$-Cartier $\Qbb$-divisor $D$ of $X$).

There exist $B$-stable rational curves $C_\mu$ and $C_{\alpha,v}$ in $X$ (which do not depend on the choice of $Q$), indexed by edges $\mu$ of $Q$ and by pairs $(\alpha,v)$ with $\alpha\in S\backslash R$ and $v$ a vertex of $Q$ that is not in the wall $W_{\alpha,P}$, such that:
\begin{enumerate}
\item the classes $[C_\mu]$ and $[C_{\alpha,v}]$ of these curves generate $NE(X)$, which is then closed and polyhedral;
\item for any edge $\mu$ of the moment polytope $Q$, the intersection number $D.C_\mu$ is the integral length of $\mu$, i.e., the length of $\mu$ divided by the length of the primitive element in the direction of $\mu$;
\item for any pair $(\alpha,v)$ as above, we have $D.C_{\alpha,v}=\langle v,\alpha^\vee\rangle$.
\end{enumerate}
\end{prop}

\section{Log MMP via a one-parameter family of polytopes}\label{sec:FamilyConstruction}
We begin with any horospherical pair $(X,\Delta)$ (see Definition~\ref{def:horopair}).

\subsection{The one-parameter family of polytopes}\label{sec:TheFamily}

We construct the one-parameter family of polytopes that permits to run the Log MMP for $(X,\Delta)$ as follows.  We use the notation of Section~\ref{sec:horopoly}.

We write $\Delta=\sum_{i=0}^r\delta_iX_i+\sum_{\alpha\in S\backslash R}\delta_\alpha D_\alpha$. Moreover, an anticanonical divisor of $X$ is $-K_X=\sum_{i=1}^rX_i+\sum_{\alpha\in S\backslash R}a_\alpha D_\alpha$ where $a_\alpha=\langle 2\rho^P,\alpha^\vee\rangle$ such that $\rho^P$ is the sum of positive roots of $G$ that are not roots of $P$ \cite{briondiv}.\\

Let $D$ be an ample $\Qbb$-Cartier $\Qbb$-divisor of $X$.\\

For any rational number $\epsilon>0$ small enough, the divisor $D+\epsilon (K_X+\Delta)$ is still ample (and $\Qbb$-Cartier by hypothesis), so that $(X,D+\epsilon (K_X+\Delta))$ defines a moment polytope $Q^\epsilon$ and a pseudo-moment polytope $\tilde{Q}^\epsilon$. Then we extend the definition for any rational number $\epsilon>0$ as follows: $\tilde{Q}^\epsilon:=\{x\in M_\Qbb\,\mid\, Ax\geq \tilde{B}+\epsilon \tilde{C}\}$ and $Q^\epsilon:= v^\epsilon+\tilde{Q}^\epsilon$ where the matrices $A$, $\tilde{B}$, $\tilde{C}$ and the vector $v^\epsilon$ are defined below.\\

Recall that $x_1,\dots,x_r$ denote the primitive elements of $N$ associated to the $G$-stable prime divisors $X_i$ of $X$.  Choosing an order in $S\backslash R$ we denote by $\alpha_1,\dots,\alpha_s$ its elements. We fix a basis $\mathcal{B}$ of $M$ and we denote by $\mathcal{B}^\vee$ the dual basis in $N$. Recall also that the rank of $M$ is denoted by~$n$.

Now define $A\in\mathcal{M}_{r+s,n}(\Qbb)$ whose first $r$ rows are the coordinates of the vectors $x_i$ in the basis $\mathcal{B}^\vee$ with $i\in\{1,\dots,r\}$ and whose last $s$ rows are the coordinates of the vectors $\alpha^\vee_{j|M}$ in $\mathcal{B}^\vee$ with $j\in\{1,\dots,s\}$. 

Let $\tilde{B}$ be the column matrix such that the pseudo-moment polytope of $D$ is defined by $\{x\in M_\Qbb\,\mid\, Ax\geq \tilde{B}\}$. In fact, if $D=\sum_{i=1}^r b_iX_i+\sum_{\alpha\in S\backslash R}b_\alpha D_\alpha$, then $\tilde{B}$ is the column matrix associated to the vector $(-b_1,\dots,-b_r,-b_{\alpha_1},\dots,-b_{\alpha_s})$. 

Similarly, the column matrix $\tilde{C}$ corresponds to the vector $(1-\delta_1,\dots,1-\delta_r,a_{\alpha_1}-\delta_{\alpha_1},\dots,a_{\alpha_s}-\delta_{\alpha_s})$.

Finally, define $v^\epsilon:=\sum_{\alpha\in S\backslash R}(b_\alpha+\epsilon (\delta_\alpha-a_\alpha))\varpi_\alpha$ (which is not necessarily in $M_\Qbb$).\\

Note that, if $M_\Qbb=X(P)_\Qbb$, we can also write the moment polytopes as follows
$Q^\epsilon:=\{x\in M_\Qbb\,\mid\, Ax\geq B+\epsilon C\}$ where $B$ and $C$ are respectively the column matrices associated to the vectors $(-b_1+\langle x_1,\sum_{\alpha\in S\backslash R}b_\alpha\varpi_\alpha\rangle,\dots,-b_r+\langle x_r,\sum_{\alpha\in S\backslash R}b_\alpha\varpi_\alpha\rangle,0,\dots,0)$ and $(1-\delta_1+\langle x_1,\sum_{\alpha\in S\backslash R}(\delta_\alpha-a_\alpha)\varpi_\alpha\rangle,\dots,1-\delta_r+\langle x_r,\sum_{\alpha\in S\backslash R}(\delta_\alpha-a_\alpha)\varpi_\alpha\rangle,0,\dots,0)$. Moreover, even if $M_\Qbb\neq X(P)_\Qbb$, it is easy to see that the $s$ last inequalities defining $\tilde{Q}^\epsilon$ are equivalent to the fact that $Q^\epsilon$ is in $X^+(P)$.

\begin{ex}\label{ex:familyII}
We consider Example~\ref{ex:GH}~{\it\large II} and the projective polarized $G/H$-embedding $(X,D)$ associated to the moment polytope $Q_7$ of Example~\ref{ex:Qi}, and the pseudo-moment polytope $-v^0+Q_7$ with $v^0=4(\varpi_\alpha+\varpi_\beta)$. Recall that $M$ is of rank $n=1$, generated by $\varpi_\alpha+2\varpi_\beta$. Choose $\mathcal{B}=(\varpi_\alpha+2\varpi_\beta)$. Then $A=\left(\begin{array}{c}1 \\-1\\1\\2
\end{array}\right)$ and $\tilde{B}=\left(\begin{array}{c}-1 \\-1\\-4\\-4
\end{array}\right)$. Moreover $-K_X=X_1+X_2+2D_\alpha+2D_\beta$.

If $\Delta=0$ we get $\tilde{C}=\left(\begin{array}{c}1 \\1\\2\\2
\end{array}\right)$ and $v^\epsilon=(4-2\epsilon)(\varpi_\alpha+\varpi_\beta)$ (which is not in $M_\Qbb$, except for $\epsilon=2$). We can easily compute that, for any $\epsilon\in[0,1]$, $\tilde{Q}^\epsilon$ is the interval $[-1+\epsilon,1-\epsilon]$ of $M_\Qbb\simeq\Qbb$, and then the $Q^\epsilon$'s are the intervals represented in the following picture, reduced to a point if $\epsilon=1$ (the intervals are more and more gray as $\epsilon$ grows).
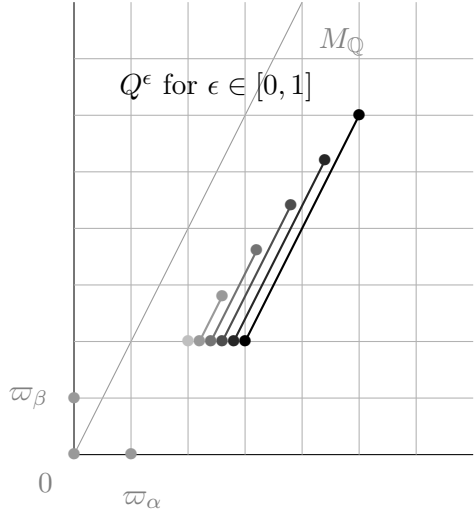
\begin{center}
\begin{tikzpicture}[scale=0.75]
\chadomII
\draw[thick] (5,6) -- (3,2);
\node at (5,6)  {\textbullet};
\node at (3,2)  {\textbullet};

\draw[thick,color=black!85] (5-0.6,6-0.8) -- (3-0.2,2);
\node at (5-0.6,6-0.8) [color=black!85] {\textbullet};
\node at (3-0.2,2) [color=black!85] {\textbullet};

\draw[thick,color=black!70] (5-1.2,6-1.6) -- (3-0.4,2);
\node at (5-1.2,6-1.6) [color=black!70] {\textbullet};
\node at (3-0.4,2) [color=black!70] {\textbullet};

\draw[thick,color=black!55] (5-1.8,6-2.4) -- (3-0.6,2);
\node at (5-1.8,6-2.4) [color=black!55] {\textbullet};
\node at (3-0.6,2) [color=black!55] {\textbullet};

\draw[thick,color=black!40] (5-2.4,6-3.2) -- (3-0.8,2);
\node at (5-2.4,6-3.2) [color=black!40] {\textbullet};
\node at (3-0.8,2) [color=black!40] {\textbullet};

\node at (2,2) [color=black!25] {\textbullet};

\node at (2.5,6.5) {$Q^\epsilon$ for $\epsilon\in[0,1]$};
\end{tikzpicture}
\end{center}

If $\Delta=X_1+X_2$ we get $\tilde{C}=\left(\begin{array}{c}0 \\0\\2\\2
\end{array}\right)$ and still $v^\epsilon=(4-2\epsilon)(\varpi_\alpha+\varpi_\beta)$. We can easily compute that, for any $\epsilon\in[0,1]$, $\tilde{Q}^\epsilon=[-1,1]$; for any $\epsilon\in[1,2]$, $\tilde{Q}^\epsilon=[-2+\epsilon,1]$;  and for any $\epsilon\in[2,\frac{5}{2}]$, $\tilde{Q}^\epsilon=[-4+2\epsilon,1]$ . Then the $Q^\epsilon$'s are the intervals represented in the following picture, reduced to the point $\varpi_\beta$ if $\epsilon=\frac{5}{2}$ (the intervals are more and more gray as $\epsilon$ grows).
\begin{center}
\begin{tikzpicture}[scale=0.75]
\chadomII
\draw[thick] (5,6) -- (3,2);
\node at (5,6)  {\textbullet};
\node at (3,2)  {\textbullet};

\draw[thick,color=black!92] (4.5,5.5) -- (2.5,1.5);
\node at (4.5,5.5) [color=black!92] {\textbullet};
\node at (2.5,1.5) [color=black!92] {\textbullet};

\draw[thick,color=black!84] (4,5) -- (2,1);
\node at (4,5) [color=black!84] {\textbullet};
\node at (2,1) [color=black!84] {\textbullet};

\draw[thick,color=black!76] (3.5,4.5) -- (1.5,0.5);
\node at (3.5,4.5) [color=black!76] {\textbullet};
\node at (1.5,0.5) [color=black!76] {\textbullet};

\draw[thick,color=black!68] (3,4) -- (1,0);
\node at (3,4) [color=black!68] {\textbullet};
\node at (1,0) [color=black!68] {\textbullet};

\draw[thick,color=black!60] (2.5,3.5) -- (0.75,0);
\node at (2.5,3.5) [color=black!60] {\textbullet};
\node at (0.75,0) [color=black!60] {\textbullet};

\draw[thick,color=black!52] (2,3) -- (0.5,0);
\node at (2,3) [color=black!52] {\textbullet};
\node at (0.5,0) [color=black!52] {\textbullet};

\draw[thick,color=black!44] (1.5,2.5) -- (0.25,0);
\node at (1.5,2.5) [color=black!44] {\textbullet};
\node at (0.25,0) [color=black!44] {\textbullet};

\draw[thick,color=black!36] (1,2) -- (0,0);
\node at (1,2) [color=black!36] {\textbullet};
\node at (0,0) [color=black!36] {\textbullet};

\draw[thick,color=black!28] (0.5,1.5) -- (0,0.5);
\node at (0.5,1.5) [color=black!28] {\textbullet};
\node at (0,0.5) [color=black!28] {\textbullet};

\node at (0,1) [color=black!20] {\textbullet};

\node at (2.5,6.5) {$Q^\epsilon$ for $\epsilon\in[0,\frac{5}{2}]$};
\end{tikzpicture}
\end{center}
\end{ex}

\begin{rem}\label{rem:deffamily}
In \cite[Section~3.1]{MMPhoro}, the definition of the family of polytopes seems to be more involved. Indeed we extended step by step the family of pseudo-moment polytopes to any rational number $\epsilon>0$, by deleting a row of $A$, $\tilde{B}$ and $\tilde{C}$ as soon as this row does not correspond to a facet of $\tilde{Q}^\epsilon$. It gives the same family of polytopes because of the convexity of the set of $\epsilon$'s such that a specific row corresponds to a facet of $\tilde{Q}^\epsilon$. But we had to give this other construction to define the good equivalence relation in this family of polytopes so that the equivalence classes of $G/H$-polytopes (Definition~\ref{def:G/H-eq}) in the family $(Q^\epsilon)_{\epsilon\in[0,\epsilon_{max}[}$ and the equivalence classes (\cite[Definition~3.14]{MMPhoro}) in the family $(\tilde{Q}^\epsilon)_{\epsilon\in[0,\epsilon_{max}[}$ are the same (\cite[Proposition~4.1]{MMPhoro}), where $\epsilon_{max}$ is the smallest (may be $+\infty$) non-negative rational number $\epsilon$ such that $Q^\epsilon$ is not a $G/H$-polytope. If $\epsilon_{max}$ is finite, it is a positive rational number by \cite[Corollary~3.16]{MMPhoro}.
 \end{rem}
 
 \begin{rem}\label{rem:endifeffective}
Note also that, if $\tilde{C}\geq 0$ and non-zero, $\epsilon_{max}$ is finite (and rational). Indeed, since $\tilde{C}\geq 0$, we get easily that for any $\epsilon>0$,  $\tilde{Q}^\epsilon\subset\tilde{Q}^0$. But $\tilde{C}\neq 0$ and $\tilde{Q}^0$ is bounded, then there exists an index $i$ in $\{1,\dots,m+r\}$ such that $\tilde{C}_i>0$, so that the set $\tilde{Q}^0\cap\{X\in M_\Qbb\,\mid\,A_iX\geq \tilde{B}_i+\epsilon\tilde{C}_i\}$ is empty for $\epsilon$ big enough (even if $A_i=0$); where $A_i$ (respectively $\tilde{B}_i$ and $\tilde{C}_i$) denotes the i-th row of $A$ (respectively of $\tilde{B}$ and $\tilde{C}$). Hence, $\tilde{Q}^\epsilon$ is empty for $\epsilon$ big enough, in particular $Q^\epsilon$ cannot be a $G/H$-polytope for any $\epsilon>0$.

Moreover, if $\epsilon_{max}$ is finite, $Q^{\epsilon_{max}}$ is neither the empty set nor a $G/H$-polytope and for any $\epsilon>\epsilon_{max}$, $Q^{\epsilon}$ is the empty set \cite[Remark~3.18]{MMPhoro}.
\end{rem}

\begin{rem}\label{rem:QCartier}
The construction above of the families $(Q^\epsilon)_{\epsilon\geq 0}$ and $(\tilde{Q}^\epsilon)_{\epsilon\geq 0}$ can be made for any $B$-stable $\Qbb$-divisor $D'$ instead of $K_X+\D$. And then, $D'$ is $\Qbb$-Cartier if and only if for $\epsilon$ small enough the polytopes $Q^\epsilon$ are $G/H$-polytopes equivalent to $Q^0$. Indeed, if $D'$ is $\Qbb$-Cartier, then for $\epsilon$ small enough, $D+\epsilon D'$ is $\Qbb$-Cartier and ample, and $Q^\epsilon$ is the moment polytope of $(X,D+\epsilon D')$, in particular it is equivalent to the moment polytope $Q^0$ of $(X,D)$. Inversely, let $\epsilon$ small enough such that $Q^\epsilon$ is a $G/H$-polytope equivalent to $Q^0$. Then the pair $(Q^\epsilon,\tilde{Q}^\epsilon)$ corresponds to a unique polarized variety $(X,D'')$ where $D''$ is an ample $\Qbb$-Cartier $B$-stable $\Qbb$-divisor. Using 4.~of Proposition~\ref{prop:moment}, we can compute that  $D''=D+\epsilon D'$, in particular $D'$ is $\Qbb$-Cartier.
\end{rem}

\subsection{Construction of pairs, contractions and flips}\label{sec:constructionvarietes}

\begin{prop}\cite[Corollary~3.16]{MMPhoro}
\label{prop:decoupage}
Let $(Q^\epsilon)_{\epsilon\in\Qbb_{\geq 0}}$ be the family of polytopes constructed above.
Then there exist non-negative integers $k,j_0,\dots,j_k$, rational numbers $\alpha_{i,j}$ for $i\in\{0,\dots,k\}$ and $j\in\{0,\dots,j_i+1\}$ and $\alpha_{k,j_k+1}=\epsilon_{max}\in\Qbb_{>0}\cup\{+\infty\}$, such that  $\alpha_{i,j_i+1}=\alpha_{i+1,0}$ for any $i\in\{0,\dots,k-1\}$, and ordered as follows,
\begin{enumerate}
\item $\alpha_{0,0}=0$;
\item for any $i\in\{0,\dots,k\}$, and for any $j<j'$ in $\{0,\dots,j_i+1\}$ we have $\alpha_{i,j}<\alpha_{i,j'}$;
\end{enumerate}

such that the (isomorphism classes of) $G/H$-embeddings that correspond to a moment polytope in the family  $(Q^\epsilon)_{\epsilon\in[0,\epsilon_{max}[}$ are:
\begin{enumerate}
\item $X_{i,j}$ for any  $i\in\{0,\dots,k\}$ and $j\in\{0,\dots,j_i\}$, respectively associated to moment polytopes $Q^\epsilon$ with $\epsilon\in[\alpha_{i,0},\alpha_{i,1}[$ if $j=0$ and $\epsilon\in]\alpha_{i,j},\alpha_{i,j+1}[$ if $j>0$;
\item $Y_{i,j}$ for any $i\in\{0,\dots,k\}$ and $j\in\{1,\dots,j_i\}$, respectively associated to the moment polytope $Q^{\alpha_{i,j}}$.
\end{enumerate}
\end{prop}

Note that $X_{0,0}$ is the variety $X$ from which we construct the family $(Q^\epsilon)_{\epsilon\in\Qbb_{\geq 0}}$.

\begin{ex}
Consider the first family of Example~\ref{ex:familyII}. Then $k=0$, $j_0=0$, $\alpha_{0,0}=0$, and $\alpha_{0,1}=1=\epsilon_{max}$. In particular, there is only one $G/H$-embedding $X_{0,0}=X$ that corresponds to a moment polytope in the family $(Q^\epsilon)_{\epsilon\in[0,1[}$, i.e., all the $G/H$-polytopes of the family $(Q^\epsilon)_{\epsilon\in[0,1[}$ are equivalent (see Definition~\ref{def:G/H-eq}) and then correspond to the same $G/H$-embedding (by Proposition~\ref{prop:G/H-eq}).

Consider now the second family of Example~\ref{ex:familyII}. Then $k=1$, $j_0=0$, $j_1=1$, $\alpha_{0,0}=0$, $\alpha_{0,1}=\alpha_{1,0}=1$, $\alpha_{1,1}=2$ and $\alpha_{1,2}=\frac{5}{2}=\epsilon_{max}$. In particular, the $G/H$-embeddings that correspond to a moment polytope in the family  $(Q^\epsilon)_{\epsilon\in[0,\frac{5}{2}[}$ are  $X_{0,0}=X$, $X_{1,0}$, $Y_{1,1}$ and $X_{1,1}$ respectively associated to the moment polytopes represented in the following figure.
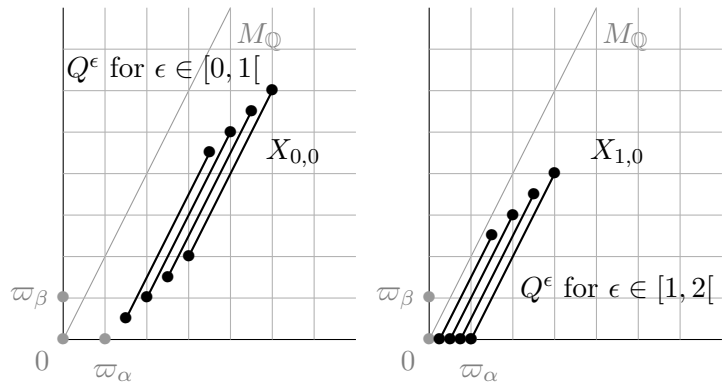
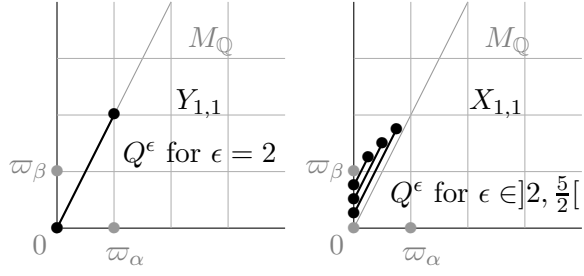
\begin{center}
\begin{tikzpicture}[scale=0.55]
\chadomII
\draw[thick] (5,6) -- (3,2);
\node at (5,6)  {\textbullet};
\node at (3,2)  {\textbullet};

\draw[thick,color=black] (4.5,5.5) -- (2.5,1.5);
\node at (4.5,5.5) [color=black] {\textbullet};
\node at (2.5,1.5) [color=black] {\textbullet};

\draw[thick,color=black] (4,5) -- (2,1);
\node at (4,5) [color=black] {\textbullet};
\node at (2,1) [color=black] {\textbullet};

\draw[thick,color=black] (3.5,4.5) -- (1.5,0.5);
\node at (3.5,4.5) [color=black] {\textbullet};
\node at (1.5,0.5) [color=black] {\textbullet};

\node at (5.5,4.5) {$X_{0,0}$};

\node at (2.4,6.5) {$Q^\epsilon$ for $\epsilon\in[0,1[$};
\end{tikzpicture}
\begin{tikzpicture}[scale=0.55]
\chadomII

\draw[thick,color=black] (3,4) -- (1,0);
\node at (3,4) [color=black] {\textbullet};
\node at (1,0) [color=black] {\textbullet};

\draw[thick,color=black] (2.5,3.5) -- (0.75,0);
\node at (2.5,3.5) [color=black] {\textbullet};
\node at (0.75,0) [color=black] {\textbullet};

\draw[thick,color=black] (2,3) -- (0.5,0);
\node at (2,3) [color=black] {\textbullet};
\node at (0.5,0) [color=black] {\textbullet};

\draw[thick,color=black] (1.5,2.5) -- (0.25,0);
\node at (1.5,2.5) [color=black] {\textbullet};
\node at (0.25,0) [color=black] {\textbullet};

\node at (4.5,4.5) {$X_{1,0}$};

\node at (4.5,1.2) {$Q^\epsilon$ for $\epsilon\in[1,2[$};
\end{tikzpicture}

\begin{tikzpicture}[scale=0.75]
\chadomIIter

\draw[thick,color=black!] (1,2) -- (0,0);
\node at (1,2) [color=black!] {\textbullet};
\node at (0,0) [color=black!] {\textbullet};

\node at (2.5,2.2) {$Y_{1,1}$};

\node at (2.5,1.3) {$Q^\epsilon$ for $\epsilon=2$};
\end{tikzpicture}
\begin{tikzpicture}[scale=0.75]
\chadomIIter

\draw[thick,color=black!] (0.75,1.75) -- (0,0.25);
\node at (0.75,1.75) [color=black!] {\textbullet};
\node at (0,0.25) [color=black!] {\textbullet};

\draw[thick,color=black!] (0.5,1.5) -- (0,0.5);
\node at (0.5,1.5) [color=black!] {\textbullet};
\node at (0,0.5) [color=black!] {\textbullet};

\draw[thick,color=black!] (0.25,1.25) -- (0,0.75);
\node at (0.25,1.25) [color=black!] {\textbullet};
\node at (0,0.75) [color=black!] {\textbullet};

\node at (2.5,2.2) {$X_{1,1}$};

\node at (2.3,0.6) {$Q^\epsilon$ for $\epsilon\in]2,\frac{5}{2}[$};
\end{tikzpicture}
\end{center}
\end{ex}

All the varieties of Proposition~\ref{prop:decoupage} have the same $B$-stable and not $G$-stable prime divisors as $X$ because they have the same open $G$-orbit $G/H$. Moreover, any of their
$G$-stable prime divisors corresponds to a $G$-stable prime divisor of $X$ (but the inverse is not necessarily true). Indeed a $G$-stable prime divisor of one of these $G/H$-embeddings corresponds to a facet $F$, of the corresponding moment polytope, that is not in any wall of the dominant chamber. In particular, if $\tilde{F}$ denotes the facet of $\tilde{Q}$ obtained by translation of $F$, $\tilde{F}$ is defined by a hyperplane associated to one of the $r$ first rows of $A$, $\tilde{B}$ and $\tilde{C}$. We still denote by $X_k$, with $k\in\{1,\dots,r\}$, the $G$-stable prime divisor of $X_{i,j}$ or $Y_{i,j}$ when it is still a divisor of $X_{i,j}$ or $Y_{i,j}$. 

Then, we define 
$$\Delta_{X_{i,j}}:=\sum_{k,\,X_k\mbox{ is a divisor of }X_{i,j}}\delta_kX_k+\sum_{\alpha\in S\backslash R}\delta_\alpha D_\alpha$$ and $\Delta_{Y_{i,j}}$ similarly.

Note that for $\Delta=K_X$ we have $\Delta_{X_{i,j}}=K_{X_{i,j}}$ and $\Delta_{Y_{i,j}}=K_{Y_{i,j}}$.

\begin{rem}
For any $i$ and $j$, for any $\epsilon\in]\alpha_{i,j},\alpha_{i,j+1}[$ the pair $(Q^\epsilon,\tilde{Q}^\epsilon)$ corresponds to the polarized variety $(X_{i,j},D+\epsilon(K_{X_{i,j}}+\Delta_{X_{i,j}}))$ by the bijection of Corollary~\ref{cor:polytopepolarized}, and  $(Q^{\alpha_{i,j}},\tilde{Q}^{\alpha_{i,j}})$ corresponds to the polarized variety $(Y_{i,j},D+\epsilon(K_{Y_{i,j}}+\Delta_{Y_{i,j}}))$.
\end{rem}

If $\epsilon_{max}$ is finite, the polytope $Q^{\epsilon_{max}}$ also defines a projective horospherical $G$-variety $Z$. Indeed, we can apply Corollary~\ref{cor:polytopepolarized} to a quadruple $(P^1,M^1,Q^{\epsilon_{max}},\tilde{Q}^{\epsilon_{max}})$ to get a polarized horospherical variety. We define $P^1$ and $M^1$ such that $(P^1,M^1,Q^{\epsilon_{max}},\tilde{Q}^{\epsilon_{max}})$ is a moment quadruple as follows: $M^1_\Qbb$ is the minimal vector subspace of $M_\Qbb$ containing $\tilde{Q}^{\epsilon_{max}}$ (i.e., $M^1_\Qbb$ is the direction of $Q^{\epsilon_{max}}$) and then $M^1:=M^1_\Qbb\cap M$; $P^1$ is the parabolic subgroup containing $B$ with simple roots $R^1$ that is the union of $R$ with the set of $\alpha\in S\backslash R$ such that $Q^{\epsilon_{max}}$ is contained in the wall $W_{\alpha,P}$. 
Then the moment quadruple $(P^1,M^1,Q^{\epsilon_{max}},\tilde{Q}^{\epsilon_{max}})$ corresponds to a polarized horospherical variety $(Z,D_Z)$. In particular, $Z$ is a $G/H^1$-embedding where $H^1$ is the subgroup of $P^1$ defined as the intersection of kernels of the characters of $P^1$ in $M^1$.

Remark that, by definition, $M^1\subset M$ and $R\subset R^1$ so that we have a projection $\pi:G/H\longrightarrow G/H^1$. Note also that the above choice of $P^1$ and $H^1$ is the unique such that $(P^1,M^1,Q^{\epsilon_{max}},\tilde{Q}^{\epsilon_{max}})$ is a moment quadruple and $\pi$ has connected fibers.\\

Now, by Proposition~\ref{prop:morph}, we get dominant $G$-equivariant morphisms: 
\begin{enumerate} \item $\phi_{i,j}:X_{i,j-1}\longrightarrow Y_{i,j}$ for any $i\in\{0,\dots,k\}$ and $j\in\{1,\dots,j_i\}$;
\item $\phi_{i,j}^+:X_{i,j}\longrightarrow Y_{i,j}$ for any $i\in\{0,\dots,k\}$ and $j\in\{1,\dots,j_i\}$; \item $\phi_i:X_{i,j_i}\longrightarrow X_{i+1,0}$ for any $i\in\{0,\dots,k-1\}$;
\item and, if $\epsilon_{max}$ is finite, $\phi:X_{k,j_k}\longrightarrow Z$.
\end{enumerate}

\begin{ex}\label{ex:MMPfamilyII}
Consider the first family of Example~\ref{ex:familyII}. Then $Q^1$ is a point in the interior of $X(P)^+$ (and $\tilde{Q}^1$ is the point 0). Then $P^1=P$, and $M^1=\{0\}$. Hence $G/H^1$ ($=Z$) is the flag variety $G/P$ (i.e., $\rm{SL}_3(\Cbb)/B$).

Consider now the second family of Example~\ref{ex:familyII}. Then $Q^\frac{5}{2}$ is the point $\varpi_\beta$, in particular it is contained in the wall $W_{\alpha,P}$. Then $P^1$ is the maximal parabolic subgroup of $G$ whose unique simple root is $\beta$, $M^1=\{0\}$ and $G/H^1$ ($=Z$) is the flag variety $G/P^1$ isomorphic to $\Pbb^2$.\\
The other morphisms $\phi_0:X_{0,0}\longrightarrow X_{1,0}$, $\phi_{1,1}:X_{1,0}\longrightarrow Y_{1,1}$ and $\phi_{1,1}^+:X_{1,1}\longrightarrow Y_{1,1}$ are morphisms that we already mentioned in Example~\ref{ex:morph}.
\end{ex}

With the same proofs as in \cite[Sections~4.3 and 4.4]{MMPhoro} (in particular by replacing the condition on $X$ to be $\Qbb$-Gorenstein by the condition on $K_X+\D$ to be $\Qbb$-Cartier), we get the following results.

\begin{prop}
\begin{enumerate} \item For any $i\in\{0,\dots,k\}$ and $j\in\{1,\dots,j_i\}$, the curves $C$ that are contracted by the morphism $\phi_{i,j}$ satisfy $(K_{X_{i,j-1}}+\Delta_{X_{i,j-1}})\cdot C<0$; for any $i\in\{0,\dots,k-1\}$, the curves $C$ that are contracted by the morphism $\phi_{i}$ satisfy $(K_{X_{i,j_i}}+\Delta_{X_{i,j_i}})\cdot C<0$; and, if $\epsilon_{max}$ is finite, the curves $C$ that are contracted by the morphism $\phi$ satisfy $(K_{X_{k,j_k}}+\Delta_{X_{k,k_j}})\cdot C<0$.
 \item For any $i\in\{0,\dots,k\}$ and $j\in\{1,\dots,j_i\}$, the curves $C$ that are contracted by the morphism $\phi_{i,j}^+$ satisfy $(K_{X_{i,j}}+\Delta_{X_{i,j}})\cdot C>0$.
 \item For any $i\in\{0,\dots,k-1\}$, the morphism $\phi_{i}$ contracts at least a $G$-stable divisor of $X_{i,j_i}$.
 \end{enumerate}
\end{prop}

\begin{prop}
The pairs $(X_{i,j},\Delta_{X_{i,j}})$ with $i\in\{0,\dots,k\}$ and $j\in\{0,\dots,j_i\}$ are horospherical pairs (i.e., $K_{X_{i,j}}+\Delta_{X_{i,j}}$ is $\Qbb$-Cartier).
(And the pairs $(Y_{i,j},\Delta_{Y_{i,j}})$  with $i\in\{0,\dots,k\}$ and $j\in\{1,\dots,j_i\}$ are not horospherical pairs (i.e., $K_{Y_{i,j}}+\Delta_{Y_{i,j}}$ is not $\Qbb$-Cartier).)
\end{prop}

Note that, for any $i\in\{0,\dots,k-1\}$, the $\Qbb$-Cartier divisor $\phi_{i}^*(K_{X_{i+1,0}}+\Delta_{X_{i+1,0}})-(K_{X_{i,j_i}}+\Delta_{X_{i,j_i}})$ is supported in the exceptional locus of $\phi_{i}$, then $\phi_{i}$ contracts a Cartier divisor.

\subsection{$\Qbb$-factorial Log MMP}

For $D$ general, the Log MMP works also for the family of $\Qbb$-factorial horospherical pairs.
In this case, all the contractions that appear above in the Log MMP are contractions of extremal rays:

\begin{prop}\label{prop:Qfactorial}
Let $(X,\Delta)$ be a horospherical pair such that $X$ is $\Qbb$-factorial.
Choose $D$ such that the vector $\tilde{B}$ is in the open set
$$\bigcup_{I\subset\{1,\dots,r+s\},\,|I|>n}\pi_I^{-1}(\Qbb^{|I|}\backslash\Im(A_I)),$$ where $A_I$ is the matrix consisting in the rows of $A$ indexed by integers in $I$ and $\pi_I$ is the canonical projection of $\Qbb^{r+s}$ to its vector subspace corresponding to the coordinates in $I$.

Then, for any $i\in\{0,\dots,k\}$ and any $j\in\{0,\dots,j_i\}$, the variety $X_{i,j}$ is $\Qbb$-factorial. 
\end{prop}

The proof is exactly the same as the proof of \cite[Proposition~4.6]{MMPhoro}.

\begin{rem}
The open set where $\tilde{B}$ is chosen, is clearly not empty and dense in $\Qbb^{r+s}$, because for any $I$ of cardinality greater than $n$, the image of $A_I$ is of codimension at least one. And, since $X$ is $\Qbb$-factorial, any vector $\tilde{B}\in\Qbb^{r+s}$ gives a $\Qbb$-Cartier divisor.
\end{rem}

\begin{prop}\label{prop:onlyrays}
Let $(X,\Delta)$ be a horospherical pair such that $X$ is $\Qbb$-factorial. 
If $D$ is general in the set of ample $\Qbb$-Cartier $\Qbb$-divisors, all morphisms $\phi_{i,j}$, $\phi_{i,j}^+$, $\phi_i$ and $\phi$ (if $\epsilon_{max}$ is finite) defined in Section~\ref{sec:constructionvarietes} are contractions of rays of the corresponding effective cones $NE(X_{i,j})$. 
\end{prop}

The proof is the same as the proof \cite[Proposition~4.8]{MMPhoro} by replacing $K_Y$ by $K_Y+\Delta_Y$ and all $K_{X_{i,j}}$ by $K_{X_{i,j}}+\Delta_{X_{i,j}}$.

\begin{rem}
Proposition~\ref{prop:onlyrays} is not true without the hypothesis of $\Qbb$-factoriality (at least for flips). In Example~\ref{ex:MMPfamilyII}, we describe a flip that consists on exchanging the color $\alpha$ for the color $\beta$. In \cite[Example~5.6]{MMPhoro} we give an example of flip (from a non-$\Qbb$-factorial variety $X$) that exchanges two colors for one color. Then the morphism $\phi^+:\,X^+\longrightarrow Y$ is obtained by adding two colors and can be factorized by one of the two morphism that adds one of the two colors. In particular $\phi^+$ is not a contraction of a ray of $NE(X^+)$ but of a 2-dimensional face of $NE(X^+)$.\\
(Such an example also exists for toric varieties of dimension~3.)
\end{rem}

\subsection{General fibers of contractions of fiber type}\label{sec:generalfiber}

We assume here that $\epsilon_{max}$ is finite.

Then, we can also describe the general fibers of the morphism $\phi:X_{k,j_{k}}\longrightarrow Z$ defined in Section~\ref{sec:constructionvarietes}. We may assume that $k=0$ and $j_0=0$, in particular $X_{k,j_k}=X$. 

\begin{teo}\label{th:generalfibers}
Let $(X,\Delta)$ be a horospherical pair, with $X$ a $G/H$-embedding. Let $D$ be an ample $\Qbb$-Cartier $B$-stable $\Qbb$-divisor on $X$. Suppose that there exists a positive rational $\epsilon_1$  such that for any $\epsilon\in[0,\epsilon_1[$ the $G/H$-polytope $Q^\epsilon$ is equivalent to $Q=Q^0$ and $Q^{\epsilon_1}$ is not a $G/H$-polytope.

Let $H^1$, $P^1$, $R^1$ and $M^1$ defined as in Section~\ref{sec:constructionvarietes}. In particular $Q^{\epsilon_1}$ is a $G/H^1$-polytope. Denote by $Z$ the associated $G/H^1$-embedding and by $\phi$ the $G$-equivariant morphism from $X$ to $Z$.

Then the general fibers of $\phi$ are either the flag variety $P^1/P$ or a projective horospherical variety $F_\phi$. Moreover in the second case, $F_\phi$ is a $L^1/H^2$-embedding, where $L^1:= H^1/R_u(H^1)$ and $H^2:=H/R_u(H^1)$ ($R_u(H^1)$ denoting the unipotent radical of $H^1$ in $G$). And a moment polytope of $F_\phi$ is the projection of $Q$ in $X(P)_\Qbb/M^1_\Qbb$.

Assume now that $X$ is $\Qbb$-factorial. Then, for general $D$, the general fibers $P^1/P$ or $F_\phi$ have Picard number one.
\end{teo}

The proof and the description of $F_\phi$ are the same as in \cite[Section~4.6]{MMPhoro}, where we only used the matrices $A$, $\tilde{B}$, $\tilde{C}$ and the vector $v^0$, and then never used that $\tilde{C}$ is defined from $K_X$, so that we can replace $K_X$ by any other $\Qbb$-Cartier divisor (under the assumptions of the beginning of the section).

\begin{ex}
Consider  the first family of Example~\ref{ex:familyII}, and the contraction of fiber type $\phi:\,X\longrightarrow G/P$. Then, since $H^1=P^1=P=B$, the group $L^1$ is isomorphic to the maximal torus $T$ of $G$, i.e., isomorphic to $(\Cbb^*)^2$. Moreover, the group $H^2$ is isomorphic to a subtorus $\Cbb^*$ of $T$. Then $L^1/H^2$ is isomorphic to $\Cbb^*$ and $F_\phi$ is isomorphic to $\Pbb^1$.

See another example in Example~\ref{ex:ex2}.
\end{ex}

\section{Conclusion and examples}\label{sec:conclusion}

Let $(X,\Delta)$ be a horospherical pair. Suppose that $-(K_X+\D)$ is non-zero and effective.
Pick an ample $\Qbb$-Cartier $B$-stable $\Qbb$-divisor $D$ of $X$. Then the family $(Q^\epsilon)_{\epsilon\geq 0}$ of polytopes defined in Section~\ref{sec:TheFamily} describes the Log MMP for $(X,\Delta)$. Moreover, it preserves some singularities of $X$ and $(X,\D)$.

First we give the definition of klt and lc singularities that we use here. It is equivalent to \cite[Definition 2.34]{KollarMori}.

\begin{defi}\label{def:kltlc}
Let $(X,\Delta)$ be a log pair such that $\lceil\D\rceil\leq 1$ (i.e., the coefficients of $\D$ are at most~1).
 A \emph{log resolution} of $(X,\Delta)$ is a proper birational map $\phi:\,V\longrightarrow X$ where $V$ is smooth and such that $\operatorname{Exc}(\phi)+\phi_*^{-1}(\D)$ is a divisor whose support has simple normal crossings.
 
 The pair  $(X,\Delta)$ has \emph{klt} (respectively \emph{lc}) singularities if there exists a log resolution $\phi:\,V\longrightarrow X$ such that every coefficient of the divisor $K_V-\phi^*(K_X+\D)$ of $V$ is strictly greater than $-1$ (respectively greater than or equal to $-1$).
\end{defi}

In the case of horospherical varieties, we can construct log resolutions by using Bott-Samelson resolutions and then we get the following characterization.
\begin{teo}\cite{KltHoro}
A horospherical pair $(X,\Delta)$ has klt singularities if and only if $\lfloor \D\rfloor\leq 0$ (i.e., the coefficient of $\D$ are less than 1).
A horospherical pair $(X,\Delta)$ has lc singularities if and only if $\lceil\D\rceil\leq 1$.
\end{teo}
In \cite{KltHoro}, we proved the first statement with the assumption that $\D$ is effective. But it is easy to see that this assumption is not necessary. And it is also not difficult to use the same proof with ``$\geq$'' instead of ``$>$'' to get the second statement.\\

 Hence, we deduce easily that the family $(Q^\epsilon)_{\epsilon\geq 0}$ gives the different types of MMP as follows.
\begin{enumerate}
\item if $X$ is $\Qbb$-factorial, we get the $\Qbb$-factorial Log MMP;
\item if $(X,\Delta)$ has klt singularities, we get the Log MMP for klt pairs (i.e., every log pair $(X_{i,j},\D_{i,j})$ has klt singularities);
\item if $(X,\Delta)$ has log canonical singularities, we get the Log MMP for lc pairs (i.e., every log pair $(X_{i,j},\D_{i,j})$ has lc singularities).
\end{enumerate}

Note also that the results above are still true if we assume in the definition of pairs that $\D$ is effective.

\begin{rem}
The assumption on $-(K_X+\D)$ to be effective is too restrictive: without this assumption, it could happen that the family $(Q^\epsilon)_{\epsilon\geq 0}$ describes the Log MMP until its end (see the end of Example~\ref{ex:ex2}). Indeed, the optimal assumption to do is that $\epsilon_{max}$ (defined in Remark~\ref{rem:deffamily}) is finite. We have  ``$-(K_X+\D)$ is non-zero and effective'' implies ``$\epsilon_{max}$ is finite'' (see Remark~\ref{rem:endifeffective}), but the inverse is not true. 

Moreover, if $\epsilon_{max}=+\infty$, the family $(Q^\epsilon)_{\epsilon\geq 0}$ describes a beginning (may be ``empty'') of the Log MMP that never ends.

Remark also that  $-(K_X+\D)$ is effective as soon as $(X,\Delta)$ has lc singularities.
\end{rem}

Now we give two examples with $G/H$ as in Example~\ref{ex:GH}~{\it I} (with corresponding notation $P$, $M$,...) and we discuss the different Log MMP we can obtain from a $G/H$-embedding. 

\begin{ex}\label{ex:ex1}
We define
$$A=\left(\begin{array}{cc}
1 & -1\\
2 & 1\\
-1 & 0\\
1 & 0
\end{array}\right)\mbox{ and }B=\left(\begin{array}{c}
-b_1\\
b_1\\
-b_3\\
0
\end{array}\right),$$ where $b_1$ and $b_3$ are rational numbers such that $b_3>0$.

We can check that the polytope $Q^0:=\{X\in M_\Qbb\,\mid\,AX\geq B\}$ is a triangle that intersects the line $\Qbb\varpi_0$ in exactly one vertex. See an example of such $Q^0$ in the figure below.

 Then the moment quadruple $(P,M,Q^0,Q^0)$ corresponds to a polarized $G/H$-embedding $(X,D)$ where the ample $\Qbb$-Cartier $B$-stable $\Qbb$-divisor $D$ is of the form $b_1X_1-b_1X_2+b_3X_3+0D_\alpha$. Note that the three $G$-stable divisors of $X$ correspond to $x_1=(1,-1)$, $x_2=(2,1)$ and $x_3=(-1,0)$ respectively, and that $\alpha^\vee_M$ is $(1,0)$, in the dual basis of the basis $(\varpi_\alpha,\varpi_0)$ of $M$.

Now define $Q^\epsilon:=\{X\in M_\Qbb\,\mid\,AX\geq B+\epsilon C\}$, where $C\in M_{4,1}(\Qbb)$. We can check that $Q^\epsilon$ is a triangle that intersects the line $\Qbb\varpi_0$ in exactly one vertex for any $\epsilon$ small enough if and only if $C_{1}+C_{2}=0$ (if not, $Q^\epsilon$ is a trapezium for any non-zero $\epsilon$ small enough). 
Here, an anticanonical divisor of $X$ is $-K_X=X_1+X_2+X_3+2D_\alpha$. In particular, we can compute with Remark~\ref{rem:QCartier} that $K_X$ is not $\Qbb$-Cartier. 

More generally, if $\D=\delta_1X_1+\delta_2X_2+\delta_3X_3+\delta_\alpha D_\alpha$ is a $B$-stable $\Qbb$-divisor of $X$, then $K_X+\D$ is $\Qbb$-Cartier if and only if $4+\delta_1+\delta_2-3\delta_\alpha=0$. In particular, if $\D$ is effective, $(X,\D)$ does never have klt or even lc singularities. 

We now consider three cases. Note that the matrix $C$ that defines the moment polytopes $Q^\epsilon$ as in Section~\ref{sec:FamilyConstruction} is $C=\left(\begin{array}{cc}
-1-\delta_1+\delta_\alpha \\
-3-\delta_2+2\delta_\alpha \\
3-\delta_3-\delta_\alpha \\
0
\end{array}\right)$.
\begin{itemize}
\item For $\D=-X_1+D_\alpha$, the pair $(X,\D)$ has lc singularities. For any $\epsilon\in[0,\frac{b_3}{2}[$, the $G/H$-polytope $Q^\epsilon$ is equivalent to $Q^0$, and $Q^\frac{b_3}{2}$ is the point $(0,b_1-\frac{b_3}{2})$.  The Log MMP gives a contraction of fiber type from $X$ to a point. Note that here $\tilde{C}
\geq 0$.

For example, if $b_1=3$ and $b_3=2$ we get the following picture.
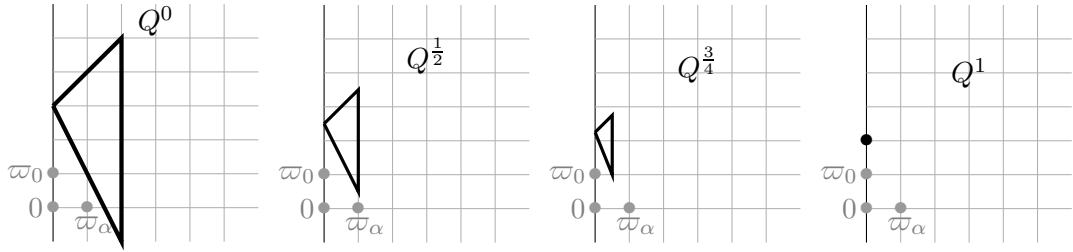
\begin{center}
\begin{tikzpicture}[scale=0.45]
\chadomIbis
\draw[ultra thick] (0,4) -- (2,6) -- (2,0) -- (0,4);
\node at (3,6.5) {$Q^0$};
\end{tikzpicture}
\begin{tikzpicture}[scale=0.45]
\chadomIbis
\draw[very thick, color=black] (0,3.5) -- (1,4.5) -- (1,1.5) -- (0,3.5);
\node at (3,5.5) [color=black] {$Q^\frac{1}{2}$};
\end{tikzpicture}
\begin{tikzpicture}[scale=0.45]
\chadomIbis
\draw[very thick, color=black] (0,3.25) -- (0.5,3.75) -- (0.5,2) -- (0,3.25);
\node at (3,5.25) [color=black] {$Q^\frac{3}{4}$};
\end{tikzpicture}
\begin{tikzpicture}[scale=0.45]
\chadomIbis

\node at (3,5) [color=black] {$Q^1$};
\node at (0,3) [color=black] {\textbullet};
\end{tikzpicture}

\end{center}

\item For $\D=X_1+X_2+2D_\alpha$, the pair $(X,\D)$ does not have lc singularities. For any $\epsilon\in[0,b_3[$, the $G/H$-polytope $Q^\epsilon$ is equivalent to $Q^0$, and $Q^{b_3}$ is the point $(0,b_1)$. The Log MMP also gives a contraction of fiber type from $X$ to a point. Note that here $\tilde{C}
\geq 0$.

For example, if $b_1=3$ and $b_3=2$ we get the following picture.
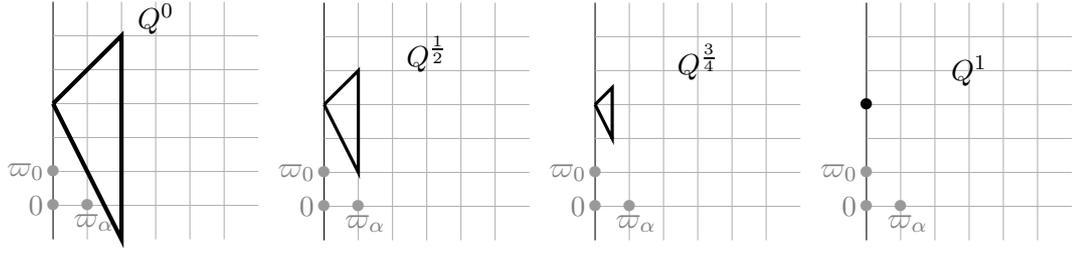
\begin{center}
\begin{tikzpicture}[scale=0.45]
\chadomIbis
\draw[ultra thick] (0,4) -- (2,6) -- (2,0) -- (0,4);
\node at (3,6.5) {$Q^0$};
\end{tikzpicture}
\begin{tikzpicture}[scale=0.45]
\chadomIbis
\draw[very thick, color=black] (0,4) -- (1,5) -- (1,2) -- (0,4);
\node at (3,5.5) [color=black] {$Q^\frac{1}{2}$};
\end{tikzpicture}
\begin{tikzpicture}[scale=0.45]
\chadomIbis
\draw[very thick, color=black] (0,4) -- (0.5,4.5) -- (0.5,3) -- (0,4);
\node at (3,5.25) [color=black] {$Q^\frac{3}{4}$};
\end{tikzpicture}
\begin{tikzpicture}[scale=0.45]
\chadomIbis

\node at (3,5) [color=black] {$Q^1$};
\node at (0,4) [color=black] {\textbullet};
\end{tikzpicture}

\end{center}

\item For $\D=\frac{5}{3}X_3+\frac{4}{3}D_\alpha$, the pair $(X,\D)$ does not have lc singularities. For any $\epsilon\in[0,+\infty[$, the $G/H$-polytope $Q^\epsilon$ is equivalent to $Q^0$ and the Log MMP does not terminate. Note that here $\tilde{C}
\ngeq 0$.

For example, if $b_1=3$ and $b_3=2$ we get the following picture.
\begin{center}
\begin{tikzpicture}[scale=0.45]
\chadomI
\draw[ultra thick] (0,4) -- (2,6) -- (2,0) -- (0,4);
\node at (3,6.5) {$Q^0$};
\end{tikzpicture}
\begin{tikzpicture}[scale=0.45]
\chadomI
\draw[very thick, color=black] (0,5) -- (2,7) -- (2,1) -- (0,5);
\node at (3,5.5) [color=black] {$Q^1$};
\end{tikzpicture}
\begin{tikzpicture}[scale=0.45]
\chadomI
\draw[very thick, color=black] (0,6) -- (2,8) -- (2,2) -- (0,6);
\node at (3,5.25) [color=black] {$Q^2$};
\end{tikzpicture}
\begin{tikzpicture}[scale=0.45]
\chadomI
\draw[very thick, color=black] (0,7) -- (2,9) -- (2,3) -- (0,7);
\node at (3,5.25) [color=black] {$Q^3$};
\end{tikzpicture}
\end{center}

\end{itemize}

\end{ex}

\begin{ex}\label{ex:ex2}
Let $A$ be as in Example~\ref{ex:ex1}. Let $B=\left(\begin{array}{c}
-b_1\\
-b_2\\
-b_3\\
0
\end{array}\right)$, where $b_1$, $b_2$ and $b_3$ are rational numbers such that $-b_1-b_2>0$ and $b_1+b_2+3b_3>0$.

We can check that the polytope $Q^0:=\{X\in M_\Qbb\,\mid\,AX\geq B\}$ is a triangle that is contained in the interior of $X(P)^+$. See an example of such $Q^0$ in the figure below.

Then the moment triple $(P,M,Q^0,Q^0)$ corresponds to the polarized $G/H$-embedding $(X,D)$ where the ample $\Qbb$-Cartier $B$-stable $\Qbb$-divisor $D$ is of the form $b_1X_1+b_2X_2+b_3X_3+0D_\alpha$. (The three $G$-stable divisors of $X$ still correspond to $x_1=(1,-1)$, $x_2=(2,1)$ and $x_3=(-1,0)$ respectively.)

Now define $Q^\epsilon:=\{X\in M_\Qbb\,\mid\,AX\geq B+\epsilon C\}$, where $C\in M_{4,1}(\Qbb)$. We can check that, for any $C$ and for any $\epsilon$ small enough, $Q^\epsilon$ is a triangle that is contained in the interior of $X(P)^+$. Hence, $X$ is $\Qbb$-factorial by Remark~\ref{rem:QCartier}. In particular, for any $B$-stable $\Qbb$-divisor $\D$ of $X$, $(X,\D)$ is a horospherical pair and $(X,0)$ has klt singularities.

One can consider the same three cases as in Example~\ref{ex:ex1} and obtain similar descriptions of the Log MMP (except that the contraction of fiber type maps to  $G/B\simeq\Pbb^1$ instead of a point).

We now consider another family of cases including the case of the pair $(X,0)$.  
Let $\D=\delta_3X_3$, with $\delta_3\in\Qbb$. We distinguish four cases (here $\tilde{C}\geq 0
$ if and only if $\delta_3\leq 1$). 
Let $\epsilon_1:=\frac{b_1+b_2+3b_3}{5-3\delta_3}$, $\epsilon_2:=\frac{-b_1-b_2}{4}$ and $\epsilon_3:=\frac{b_3}{3-\delta_3}$.
\begin{itemize}
\item If $\delta_3<\frac{5}{3}$ and $\epsilon_1<\epsilon_2$, then for any $\epsilon\in[0,\epsilon_1[$, the $G/H$-polytope $Q^\epsilon$ is equivalent to $Q^0$, and $Q^{\epsilon_1}$ is a point in the interior of the dominant chamber. Hence the Log MMP gives a contraction of fiber type to $G/B\simeq\Pbb^1$.
\item If $\delta_3<\frac{5}{3}$ and $\epsilon_1=\epsilon_2$, then for any $\epsilon\in[0,\epsilon_1[$, the $G/H$-polytope $Q^\epsilon$ is equivalent to $Q^0$, and $Q^{\epsilon_1}$ is a point in the line $\Qbb\varpi_0$. Hence the Log MMP gives a contraction of fiber type to a point.
\item If $\delta_3<\frac{5}{3}$ and $\epsilon_1>\epsilon_2$ or if $\frac{5}{3}\leq\delta_3<3$, then for any $\epsilon\in[0,\epsilon_2[$ (with $\epsilon_2$ defined above), the $G/H$-polytope $Q^\epsilon$ is equivalent to $Q^0$, $Q^{\epsilon_2}$ is a triangle that intersects the line $\Qbb\varpi_0$ in exactly one vertex, and for any $\epsilon\in]\epsilon_2,\epsilon_3[$, the $G/H$-polytope $Q^\epsilon$ is a quadrilateral that intersects the line $\Qbb\varpi_0$ along an edge, and $Q^{\epsilon_3}$ is a segment in the line $\Qbb\varpi_0$. Hence the Log MMP begins here with a flip and ends with a contraction of fiber type to the $\Cbb^*$-embedding $\Pbb^1$ ($P^1=G$ and $M^1=\Zbb\varpi_0$). The general fiber $F_\phi$ of this latter contraction is ismorphic to $\Pbb^2$: indeed, $L^1\simeq\rm{SL}_2$ and $H^2$ is isomorphic to the subgroup $U$ of $\rm{SL}_2$ consisting of upper triangular matrices with ones on the diagonal; and then $F_\phi$ is the unique projective $\rm{SL}_2/U$-embedding with the color (see \cite[Exemples~1.13~(1)]{these}).

In particular, the family $(Q^\epsilon)_{\epsilon\geq 0}$ describes the Log MMP until its end, even if $\tilde{C}$ is not non-negative, i.e., even if $-(K_X+\D)$ is not effective.  

For example, if $b_1=3$, $b_2=-4$, $b_3=2$ and $\delta_3=0$ we get the following picture. In that case $1=\epsilon_1>\epsilon_2=\frac{1}{4}$ and $\epsilon_3=\frac{2}{3}$.
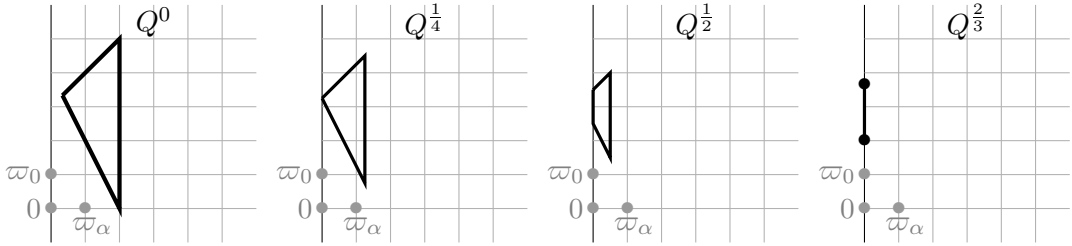
\begin{center}
\begin{tikzpicture}[scale=0.45]
\chadomIbis
\draw[ultra thick] (1/3,13/3) -- (2,6) -- (2,1) -- (1/3,13/3);
\node at (3,6.5) {$Q^0$};
\end{tikzpicture}
\begin{tikzpicture}[scale=0.45]
\chadomIbis
\draw[very thick] (1/3-4/3*1/4,13/3-1/3*1/4) -- (2-3*1/4,6-2*1/4) -- (2-3*1/4,1+3*1/4) -- (1/3-4/3*1/4,13/3-1/3*1/4);
\node at (3,6.5) {$Q^\frac{1}{4}$};
\end{tikzpicture}
\begin{tikzpicture}[scale=0.45]
\chadomIbis
\draw[very thick] (0,4+1/2) -- (2-3*1/2,6-2*1/2) -- (2-3*1/2,1+3*1/2) -- (0,5-3*1/2) -- (0,4+1/2);
\node at (3,6.5) {$Q^\frac{1}{2}$};
\end{tikzpicture}
\begin{tikzpicture}[scale=0.45]
\chadomIbis
\draw[very thick] (0,4+2/3) -- (0,5-3*2/3);
\node at (0,4+2/3) {\textbullet};
\node at (0,5-3*2/3) {\textbullet};
\node at (3,6.5) {$Q^\frac{2}{3}$};
\end{tikzpicture}
\end{center}

\item If $\delta_3\geq 3$, for any $\epsilon\in[0,\epsilon_2[$, the $G/H$-polytope $Q^\epsilon$ is equivalent to $Q^0$, $Q^{\epsilon_2}$ is a triangle that intersects the line $\Qbb\varpi_0$ in exactly one vertex, and for any $\epsilon\in]\epsilon_2,+\infty[$, the $G/H$-polytope $Q^\epsilon$ is a quadrilateral that intersects the line $\Qbb\varpi_0$ along an edge. Hence the Log MMP begins here with a flip but does not terminate.
For example, if $b_1=3$, $b_2=-4$, $b_3=2$ and $\delta_3=4$ we get the following picture. In that case $\epsilon_2=\frac{1}{4}$.
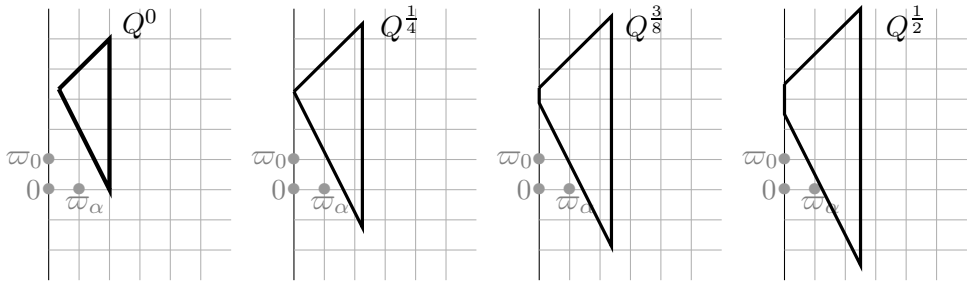
\begin{center}
\begin{tikzpicture}[scale=0.4]
\chadomIter
\draw[ultra thick] (1/3,13/3) -- (2,6) -- (2,1) -- (1/3,13/3);
\node at (3,6.5) {$Q^0$};
\end{tikzpicture}
\begin{tikzpicture}[scale=0.4]
\chadomIter
\draw[very thick] (1/3-4/3*1/4,13/3-1/3*1/4) -- (2+1/4,6+2*1/4) -- (2+1/4,1-5*1/4) -- (1/3-4/3*1/4,13/3-1/3*1/4);
\node at (3.5,6.5) {$Q^\frac{1}{4}$};
\end{tikzpicture}
\begin{tikzpicture}[scale=0.4]
\chadomIter
\draw[very thick] (0,4+3/8) -- (2+3/8,6+2*3/8) -- (2+3/8,1-5*3/8) -- (0,5-3*3/8) -- (0,4+3/8);
\node at (3.5,6.5) {$Q^\frac{3}{8}$};
\end{tikzpicture}
\begin{tikzpicture}[scale=0.4]
\chadomIter
\draw[very thick] (0,4+1/2) -- (2+1/2,6+2*1/2) -- (2+1/2,1-5*1/2) -- (0,5-3*1/2) -- (0,4+1/2);
\node at (4,6.5) {$Q^\frac{1}{2}$};
\end{tikzpicture}

\end{center}

\end{itemize} 
\end{ex}

\begin{rem}
We can imagine a Log MMP avoiding flips. For example, consider the pair $(X,0)$ of Example~\ref{ex:ex2}, and the case where $\epsilon_1>\epsilon_2$. Instead of considering the flip (i.e., $G/H$-polytopes $Q^\epsilon$ with $\epsilon\in]\epsilon_2,\frac{b_3}{3}[$), we can apply the program from the beginning to the $G/H$-embedding corresponding to the $G/H$-polytope $Q^{\epsilon_2}$, which is the $G/H$-embedding of Example~\ref{ex:ex1}.

Remark that $(X,0)$ is a pair with klt singularities with an effective $\Qbb$-divisor, but the $G/H$-embedding of Example~\ref{ex:ex1} admits no pair $(X,\D)$ with klt singularities such that $\D$ is effective.\\
But for any horospherical projective variety $X$, there exists a horospherical pair $(X,\D)$ with klt singularities. Indeed, pick any ample $\Qbb$-Cartier divisor $D'$ of $X$. Up to linear equivalence, we can suppose that $D'$ is $B$-stable and strictly effective (i.e., effective and the support of $D'$ is the union of all prime $B$-stable divisors of $X$). Then there exists a positive integer $m$ such that $\lfloor -K_X-mD'\rfloor\leq 0$, and $\D:= -K_X-mD'$ suits. Note that the constructed pair is quite special because $-(K_X+\D)=mD'$ is ample.
\end{rem}

\bibliographystyle{amsalpha}
\bibliography{logMMPviapoly}

\end{document}